\documentclass[runningheads,a4paper,abstract=on]{scrartcl}

\usepackage[utf8]{inputenc}
\usepackage[T1]{fontenc}
\usepackage{amsmath,amssymb,amsthm}
\usepackage{xspace}

\usepackage[usenames,dvipsnames,table]{xcolor}
\usepackage{graphicx}
\DeclareGraphicsExtensions{.pdf,.png,.eps}
\graphicspath{{figs/}}

\usepackage{tikz}
\usetikzlibrary{svg.path,calc}
\definecolor{orcidlogocol}{HTML}{A6CE39}
\tikzset{orcid/.pic={\node[anchor=south west] {\pgftext{\includegraphics{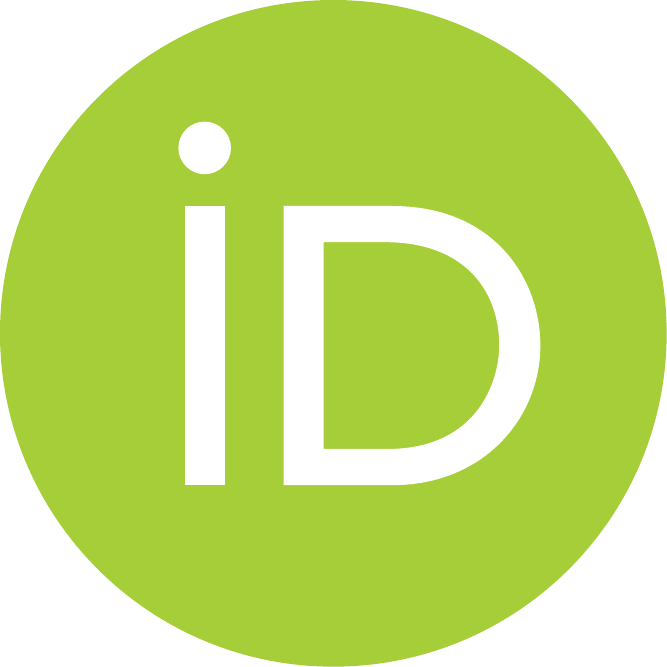}}};} x}
\newcommand\orcidicon[1]{\href{https://orcid.org/#1}{\makebox{
    \begin{tikzpicture}[ overlay,remember picture]
        \coordinate (A);
        \coordinate(B) at ($(A)-(0pt,-6pt)$);
    \end{tikzpicture}
    \begin{tikzpicture}[overlay,remember picture,yscale=0.08,xscale=0.08,transform shape]
        \pic at (B) {orcid};
    \end{tikzpicture}
}{}}}

\usepackage{multirow}
\usepackage{multicol}

\usepackage{enumitem}

\usepackage{url}
\usepackage{cite}
\usepackage[raiselinks=true,
			bookmarks=true,
			bookmarksopen=true,
			bookmarksopenlevel=2,
			bookmarksnumbered=true,
			hyperindex=true,
			plainpages=false,
			pdfpagelabels=true,
			pdfborder={0 0 0.5},
			colorlinks=false,						
			linkbordercolor={0 0.61 0.50},   
			citebordercolor={0 0.61 0.50}]{hyperref} 
\usepackage[capitalise,noabbrev,nameinlink ]{cleveref}
\newcommand{\etal}{{et~al.}}
\newcommand{\doi}[1]{\href{https://dx.doi.org/#1}{\texttt{doi:#1}}}
\newcommand{\arxiv}[1]{\href{https://arxiv.org/abs/#1}{\texttt{arXiv:#1}}}

\usepackage{authblk}
\title{Rearrangement operations on unrooted phylogenetic networks}
\date{}
\author[1]{Remie Janssen \protect\orcidicon{0000-0002-5192-1470}}
\author[2]{Jonathan Klawitter \protect\orcidicon{0000-0001-8917-5269}}
\affil[1]{Delft Institute of Applied Mathematics, Delft University of Technology, Netherlands}
\affil[2]{School of Computer Science, University of Auckland, New Zealand}

\urldef{\mailsRJ}\path|remiejanssen@gmail.com|
\urldef{\mailsJK}\path|jo.klawitter@gmail.com|
\affil[ ]{\textit {\mailsRJ}, {\mailsJK}}

\usepackage{mathtools}
\DeclarePairedDelimiter\abs{\lvert}{\rvert}
\DeclarePairedDelimiter\set{\{}{\}}
\DeclarePairedDelimiter\ceil{\lceil}{\rceil}
\DeclarePairedDelimiter\floor{\lfloor}{\rfloor}

\newtheoremstyle{break}
  {\topsep} 
  {\topsep} 
  {\itshape} 
  {0pt} 
  {\bfseries} 
  {.} 
  {\newline} 
  {} %
\theoremstyle{break} 
\newtheorem{theorem}{Theorem}[section]
\newtheorem{lemma}[theorem]{Lemma}
\newtheorem{proposition}[theorem]{Proposition}
\newtheorem{corollary}[theorem]{Corollary}
\newtheorem{observation}[theorem]{Observation}

\def\Oh{\ensuremath{\mathcal{O}}}

\newcommand{\NNI}{\textup{NNI}\xspace}
\newcommand{\NNIZ}{\textup{NNI$^0$}\xspace}
\newcommand{\NNIZimproper}{\textup{NNI$^0_{\text{improper}}$}\xspace}
\newcommand{\NNIP}{\textup{NNI$^+$}\xspace}
\newcommand{\NNIM}{\textup{NNI$^-$}\xspace}
\DeclareMathOperator{\dNNI}{\ensuremath{d_{\NNI}}}
\newcommand{\PR}{\textup{PR}\xspace}
\newcommand{\PRZ}{\textup{PR$^0$}\xspace}
\newcommand{\PRP}{\textup{PR$^+$}\xspace}
\newcommand{\PRM}{\textup{PR$^-$}\xspace}
\DeclareMathOperator{\dPR}{\ensuremath{d_{\PR}}}
\newcommand{\SNPR}{\textup{SNPR}\xspace}

\newcommand{\SPR}{\textup{SPR}\xspace}
\newcommand{\TBR}{\textup{TBR}\xspace}
\newcommand{\TBRZ}{\textup{TBR$^0$}\xspace}
\newcommand{\TBRP}{\textup{TBR$^+$}\xspace}
\newcommand{\TBRM}{\textup{TBR$^-$}\xspace}
\DeclareMathOperator{\dTBR}{\ensuremath{d_{\TBR}}}
\DeclareMathOperator{\dist}{\ensuremath{d}}
\def\distT{\ensuremath{\dist_\mathcal{T}}}
\def\distN{\ensuremath{\dist_\mathcal{N}}}
\def\utrees{\ensuremath{u\mathcal{T}_n}}
\def\utreesx[#1]{\ensuremath{u\mathcal{T}_{#1}}}
\def\unets{\ensuremath{u\mathcal{N}_n}}
\def\unetsr{\ensuremath{u\mathcal{N}_{n,r}}}
\def\unetsx[#1]{\ensuremath{u\mathcal{N}_{#1}}}
\def\unetsxx[#1,#2]{\ensuremath{u\mathcal{N}_{#1,#2}}}
\def\utbasednets{\ensuremath{u\mathcal{TB}_n}}
\def\utbasednetsr{\ensuremath{u\mathcal{TB}_{n,r}}}
\def\ulvlknets{\ensuremath{u\mathcal{LV}\text{-}k_{n}}}
\def\class{\ensuremath{\mathcal{C}_n}}

\usepackage{hyperref}
\begin{document}

\maketitle

\pdfbookmark[1]{Abstract}{Abstract} 
\begin{abstract}
Rearrangement operations transform a phylogenetic tree into another one and hence induce a metric on the space of phylogenetic trees. Popular operations for unrooted phylogenetic trees are NNI (nearest neighbour interchange), SPR (subtree prune and regraft), and TBR (tree bisection and reconnection). Recently, these operations have been extended to unrooted phylogenetic networks---generalisations of phylogenetic trees that can model reticulated evolutionary relationships---where they are called NNI, PR, and TBR moves. Here, we study global and local properties of spaces of phylogenetic networks under these three operations. In particular, we prove connectedness and asymptotic bounds on the diameters of spaces of different classes of phylogenetic networks, including tree-based and level-$k$ networks. We also examine the behaviour of shortest TBR-sequence between two phylogenetic networks in a class, and whether the TBR-distance changes if intermediate networks from other classes are allowed: for example, the space of phylogenetic trees is an isometric subgraph of the space of phylogenetic networks under TBR. Lastly, we show that computing the TBR-distance and the PR-distance of two phylogenetic networks is NP-hard. 
\end{abstract}

\section{Introduction} 
Phylogenetic trees and networks are leaf-labelled graphs that are used to visualise and study the evolutionary history of taxa like species, genes, or languages. While phylogenetic trees are used to model tree-like evolutionary histories, the more general phylogenetic networks can be used for taxa whose past includes reticulate events like hybridisation or horizontal gene transfer~\cite{SS03,HRS10,Ste16}. Such reticulate events arise in all domains of life~\cite{TN05,RW07,MMM17,WKH17}. 
In some cases, it can be useful to distinguish between rooted and unrooted phylogenetic networks. In a rooted phylogenetic network, the edges are directed from a designated root towards the leaves. Hence, it models evolution along the passing of time. An unrooted phylogenetic network, on the other hand, has undirected edges and thus represent evolutionary relatedness of the taxa. In some cases, unrooted phylogenetic networks can be thought of as rooted phylogenetic networks in which the orientation of the edges has been disregarded.
Such unrooted phylogenetic networks are called proper~\cite{JJEvIS17,FHM18}.
Here we focus on unrooted, binary, proper phylogenetic networks, where binary means that all vertices except for the leaves have degree three. The set of phylogenetic networks on the same taxa can be partitioned into tiers that contain all networks of the same size.

\begin{figure}[htb]
  \centering
  \includegraphics{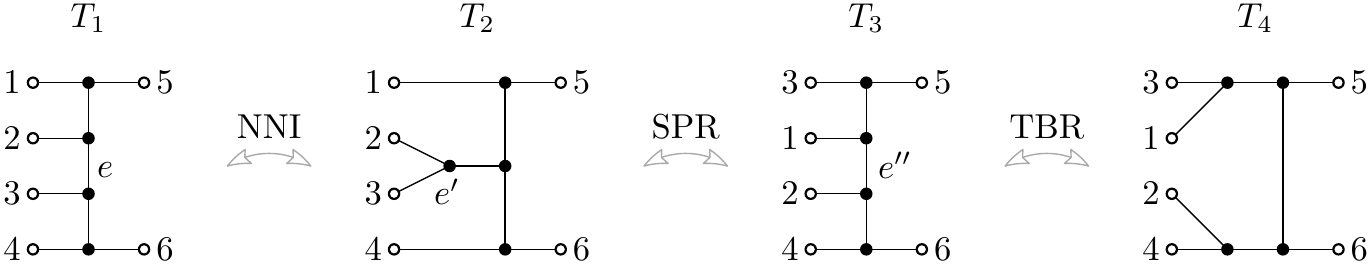}
  \caption{The three rearrangement operations on unrooted phylogenetic trees: The NNI from $T_1$ to $T_2$ changes the order of the four edges incident to $e$; the SPR from $T_2$ to $T_3$ prunes the edge $e'$, and then regrafts it again; and the TBR from $T_3$ to $T_4$ first removes the edge $e''$, and then reconnects the resulting two trees with a new edge. Note that every NNI is also an SPR and every SPR is also a TBR but not vice versa.} 
  \label{fig:utrees:rearrangementOps}
\end{figure}

A rearrangement operation transforms a phylogenetic tree into another tree by making a small graph theoretical change. An operation that works locally within the tree is the NNI (nearest neighbour interchange) operation, which changes the order of the four edges incident to an edge $e$. See for example the NNI from $T_1$ to $T_2$ in \cref{fig:utrees:rearrangementOps}.
Two further popular rearrangement operations are the SPR (subtree prune and regraft) operation, which as the name suggests prunes (cuts) an edge and then regrafts (attaches) the resulting half edge again, and the TBR (tree bisection and reconnection) operation, which first removes an edge and then adds a new one to reconnect the resulting two smaller trees. See, for example, the SPR from $T_2$ to $T_3$ and the TBR from $T_3$ to $T_4$ in \cref{fig:utrees:rearrangementOps}.

The set of phylogenetic trees on a fixed set of taxa together with a rearrangement operation yields a graph where the vertices are the trees and two trees are adjacent if they can be transformed into each other with the operation. We call this a space of phylogenetic trees. This construction also induces a metric on phylogenetic trees as the distance of two trees is then given as the distance in this space, that is, the minimum number of applications of the operation that are necessary to transform one tree into the other~\cite{SOW96}. However, computing the distance of two trees under \NNI, \SPR, and \TBR is NP-hard~\cite{DGHJLTZ97,HDRB08,AS01}. Nevertheless, both the space of phylogenetic trees and a metric on them are of importance for the many inference methods for phylogenetic trees that rely on local search strategies~\cite{Gus14,StJ17}.

Recently, these rearrangement operations have been generalised to phylogenetic networks, both for unrooted networks \cite{HLMW16,HMW16,FHMW17} and for rooted networks\cite{BLS17,FHMW17,GvIJLPS17,Kla19}.
For unrooted networks, Huber \etal~\cite{HLMW16} first generalised NNI to level-1 networks, which are phylogenetic networks where all cycles are vertex disjoint. This generalisation includes a horizontal move that changes the topology of the network, like an NNI on a tree, and vertical moves that add or remove a triangle to change the size of the network. Among other results, they then showed that the space of level-1 networks and its tiers are connected under NNI~\cite[Theorem 2]{HLMW16}. Note that connectedness implies that the distance between any two networks in such a space is finite and that NNI thus induces a metric. This NNI operation was then extended by Huber \etal~\cite{HMW16} to work for general unrooted phylogenetic networks. Again, connectedness of the space was proven. Later, Francis \etal~\cite{FHMW17} gave lower and upper bounds on the diameter (the maximum distance) of the space of unrooted phylogenetic network of a fixed size under \NNI. 
They also showed that \SPR and \TBR can straightforwardly be generalised to phylogenetic networks, that the connectedness under \NNI implies connectedness under \SPR and \TBR, and they gave bounds on the diameters. These bounds for \SPR were made asymptotically tight by Janssen \etal~\cite{JJEvIS17}. Here, we improve these bounds on the diameter under \TBR.

There are several generalisations of \SPR on rooted phylogenetic trees to rooted phylogenetic networks for which connectedness and diameters have been obtained~\cite{BLS17,FHMW17,GvIJLPS17,JJEvIS17,Jan18}. 
For example, Bordewich \etal~\cite{BLS17} introduced SNPR (subnet prune and regraft), a generalisation of \SPR that includes vertical moves, which add or remove an edge. They then proved connectedness under \SNPR for the space of rooted phylogenetic networks and for special classes of phylogenetic networks including tree-based networks. Roughly speaking, these are networks that have a spanning tree that is the subdivision of a phylogenetic tree on the same taxa~\cite{FS15,FHM18}. Furthermore, Bordewich \etal~\cite{BLS17} gave several bounds on the SNPR-distance of two phylogenetic networks. Further bounds and a characterisation of the SNPR-distance of a tree and a network were recently proven by Klawitter and Linz~\cite{KL19}. Here, we show that these bounds and characterisation on the SNPR-distance of rooted phylogenetic networks are analogous to the TBR-distance of two unrooted phylogenetic networks.

In this paper, we study spaces of unrooted phylogenetic networks under \NNI, \PR (prune and regraft), and \TBR. Here, the \PR and the \TBR operation are the generalisation of \SPR and \TBR on trees, respectively, where vertical moves add or remove an edge like the vertical moves of the SNPR operation in the rooted case. After the preliminary section, we examine the relation of \NNI, \PR, and \TBR; in particular, how a sequence using one of these operations can be transformed into a sequence using another operation (\cref{sec:relations}). We then study properties of shortest paths under \TBR in \cref{sec:paths}.
This includes the translation of the results from Bordewich \etal~\cite{BLS17} and Klawitter and Linz~\cite{KL19} on the \SNPR-distance of rooted phylogenetic networks to the \TBR-distance of unrooted phylogenetic networks.
Next, we consider the connectedness and diameters of spaces of phylogenetic networks for different classes of phylogenetic networks, including tree-based networks and level-$k$ networks (\cref{sec:connectedness}).
A subspace of phylogenetic networks (e.g., the space of tree-based networks) is an isometric subgraph of a larger space of phylogenetic networks if, roughly speaking, the distance of two networks is the same in the smaller and the larger space. In \cref{sec:isometric} we study such isometric relations and answer a question by Francis \etal~\cite{FHMW17} by showing that the space of phylogenetic trees is an isometric subgraph of the space of phylogenetic networks under \TBR. We use this result in \cref{sec:complexity} to show that computing the \TBR-distance is NP-hard. In the same section, we also show that computing the \PR-distance is NP-hard.

\section{Preliminaries} 
This section provides notation and terminology used in the remainder of the paper.
In particular, we define phylogenetic networks and special classes thereof, and rearrangement operations and how they induce distances. Throughout this paper, $X=\set{1, 2,\ldots, n}$ denotes a finite set of taxa.

\paragraph{Phylogenetic networks.}
An \emph{unrooted, binary phylogenetic network} $N$ on a set of \emph{taxa} $X$ is an undirected multigraph such that the leaves are bijectively labelled with $X$ and all non-leaf vertices have degree three. It is called \emph{proper} if every cut-edge separates two labelled leaves~\cite{FHM18}, and \emph{improper} otherwise. This property implies that every edge lies on a path that connects two leaves. More importantly, a network can be rooted at any leaf if and only if it is proper~\cite[Lemma 4.13]{JJEvIS17}. If not mentioned otherwise, we assume that a phylogenetic network is proper. Furthermore, note that our definition of a phylogenetic network permits the existence of parallel edges in $N$, i.e., we allow that two distinct edges join the same pair of vertices. An \emph{unrooted, binary phylogenetic tree} $T$ on $X$ is an unrooted, binary phylogenetic network on $X$ that is a tree. 

\begin{figure}[htb]
\centering
  \includegraphics{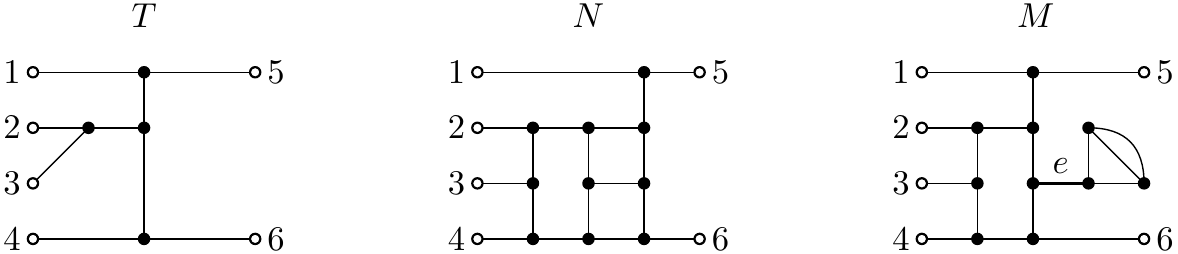}
  \caption{An unrooted, binary phylogenetic tree $T \in \utreesx[6]$ and an unrooted, binary proper phylogenetic network $N \in \unetsx[6]$. The unrooted, binary phylogenetic network $M$ is improper since the cut-edge $e$ does not lie on a path that connects two leaves.}
  \label{fig:unets:treeAndNetwork}
\end{figure}

Let $\unets$ denote the set of all unrooted, binary proper phylogenetic networks on $X$
and let $\utrees$ denote the set of all unrooted, binary phylogenetic trees on $X$, where $X = \set{1, 2,\ldots, n}$.
To ease reading, we refer to an unrooted, binary proper phylogenetic network (resp. unrooted, binary phylogenetic tree) on $X$ simply as phylogenetic network or network (resp. phylogenetic tree or tree). 
\Cref{fig:unets:treeAndNetwork} shows an example of a tree $T \in \utreesx[6]$, a network in $N \in \unetsx[6]$, 
and an improper network $M$.

An edge of a network $N$ is an \emph{external} edge if it is incident to a leaf, and an \emph{internal} edge otherwise.
A \emph{cherry} $\set{a, b}$ of $N$ is a pair of leaves $a$ and $b$ in $N$ that are adjacent to the same vertex. 
For example, each network in \cref{fig:unets:treeAndNetwork} contains the cherry $\set{1, 5}$.

\paragraph{Tiers.} 
We say a network $N = (V, E)$ has \emph{reticulation number\footnotemark} $r$ for $r = \abs{E} - (\abs{V} - 1)$,
that is, the number of edges that have to be deleted from $N$ to obtain a spanning tree of $N$.
For example, the network $N$ in \cref{fig:unets:treeAndNetwork} has reticulation number three.
Note that a phylogenetic tree is a phylogenetic network with reticulation number zero.
Let $\unetsr$ denote \emph{tier} $r$ of $\unets$, the set of networks in $\unets$ that have reticulation number $r$.
\footnotetext{In graph theory the value $\abs{E} - (\abs{V} - 1)$ of a connected graph is also called the cyclomatic number of the graph~\cite{Diestel}.}

\paragraph{Embedding.}
Let $G$ be an undirected graph. 
\emph{Subdividing} an edge $\set{u, v}$ of $G$ consists of replacing $\set{u, v}$ by a path form $u$ to $v$ that contains at least one edge. A \emph{subdivision} $G^*$ of $G$ is a graph that can be obtained from $G$ by subdividing edges of $G$. If $G$ has no degree two vertices, there exists a canonical embedding of vertices of $G$ to vertices of $G^*$ and of edges of $G$ to paths of $G^*$.
Let $N \in \unets$. We say $G$ has an \emph{embedding} into $N$ if there exists a subdivision $G^*$ of $G$ that is a subgraph of $N$ such that the embedding maps each labelled vertex of $G^*$ to a labelled vertex of $N$ with the same label. 

\paragraph{Displaying.}
Let $T \in \utrees$ and $N \in \unets$.
We say $N$ \emph{displays} $T$ if $T$ has an embedding into $N$.
For example, in \cref{fig:unets:treeAndNetwork} the tree $T$ is displayed by both networks $N$ and $M$.
Let $D(N)$ be the set of trees in $\utrees$ that are displayed by $N$.
This notion can be extended to trees with fewer leaves, and to networks.
For this, let $M$ be a phylogenetic network on $Y \subseteq X = \set{1, \ldots, n}$.
We say $N$ \emph{displays} $M$ if $M$ has an embedding into $N$.
Let $P = \set{M_1,\ldots, M_k}$ be a set of phylogenetic networks $M_i$ on $Y_i \subseteq X = \set{1, \ldots, n}$.
Then let $\unets(P)$ denote the subset of networks in $\unets$ that display each network in $P$.

\paragraph{Tree-based networks.}
A phylogenetic network $N \in \unets$ is a \emph{tree-based} network if there is a tree $T \in \utrees$ that has an embedding into $N$ as a spanning tree. In other words, there exists a subdivision $T^*$ of $T$ that is a spanning tree of $N$.
The tree $T$ is then called a \emph{base tree} of $N$. Let $\utbasednets$ denote the set of tree-based networks in $\unets$.
For $T \in \utrees$, let $\utbasednets(T)$ denote the set of tree-based networks in $\utbasednets$ with base tree $T$.

\paragraph{Level-$k$ networks.}
A blob $B$ of a network $N \in \unets$ is a nontrivial two-connected component of $N$. The \emph{level} of $B$ is the minimum number of edges that have to be removed from $B$ to make it acyclic. The \emph{level} of $N$ is the maximum level of all blobs of $N$. If the level of $N$ is at most $k$, then $N$ is called a \emph{level-$k$} network. Let $\ulvlknets$ denote the set of level-$k$ networks in $\unets$.

\begin{figure}[htb]
\centering
  \includegraphics{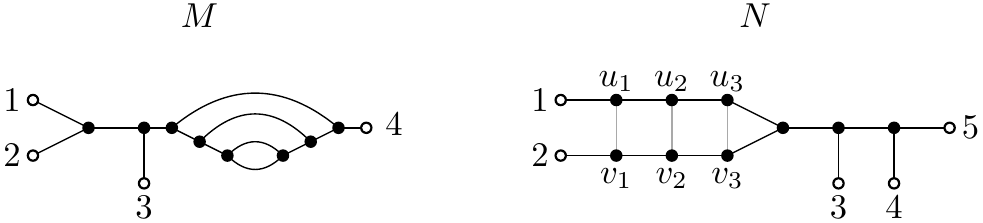}
  \caption{A network $M$ with a $3$-burl and a sorted $3$-handcuffed caterpillar $N$.}
  \label{fig:unets:handcuffed}
\end{figure}

\paragraph{$r$-Burl.} An $r$-burl is a specific type of blob that we define recursively: a $1$-burl is the blob consisting of a pair of parallel edges; an $r$-burl is the blob obtained by placing a pair of parallel edges on one of the parallel edges of an $r-1$-burl for all $r>1$. See for example the network $M$ in \cref{fig:unets:handcuffed}.

\paragraph{$r$-Handcuffed trees and caterpillars.}
Let $T \in \unets$ and let $a$ and $b$ be two leaves of $T$. Let $e$ and $f$ be the edges incident to $a$ and $b$, respectively. Subdivide $e$ and $f$ with vertices $\set{u_1, \ldots, u_r}$ and $\set{v_1, \ldots, v_r}$, respectively, and add the edges $\set{u_1, v_1}, \ldots, \set{u_r, v_r}$. The resulting network is an \emph{$r$-handcuffed tree} $N \in \unets$ with base tree $T$ on the \emph{handcuffed} leaves $\set{a, b}$. Note that $N$ has reticulation number $r$. If the tree $T$ is a caterpillar and $a$ and $b$ form a cherry of $T$, then the resulting network $N$ is an \emph{$r$-handcuffed caterpillar}. Furthermore, we call an $r$-handcuffed caterpillar \emph{sorted} if it is handcuffed on the leafs 1 and 2 and the leafs from 3 to $n$ have a non-decreasing distance to leaf 1. See \cref{fig:unets:handcuffed} for an example.

\paragraph{Suboperations.}
To define rearrangement operations on phylogenetic networks, we first define several suboperations. Let $G$ be an undirected graph. A degree-two vertex $v$ of $G$ with adjacent vertices $u$ and $w$ gets \emph{suppressed} by deleting $v$ and its incident edges, and adding the edge $\set{u, w}$. The reverse of this suppression is the subdivision of $\set{u, w}$ with vertex $v$.

Let $N \in \unets$ be a network, and $\set{u, v}$ an edge of $N$. Then $\set{u, v}$ gets \emph{removed} by deleting $\set{u, v}$ from $N$ and suppressing any resulting degree-two vertices. We say $\set{u, v}$ gets \emph{pruned} at $u$ by transforming it into the half edge $\set{\cdot, v}$ and suppressing $u$ if it becomes a degree-two vertex. Note that otherwise $u$ is a leaf. In reverse, we say that a half edge $\set{\cdot, v}$ gets \emph{regrafted} to an edge $\set{x, y}$ by transforming it into the edge $\set{u, v}$ where $u$ is a new vertex subdividing $\set{x, y}$. 

\begin{figure}[htb]
  \centering
  \includegraphics{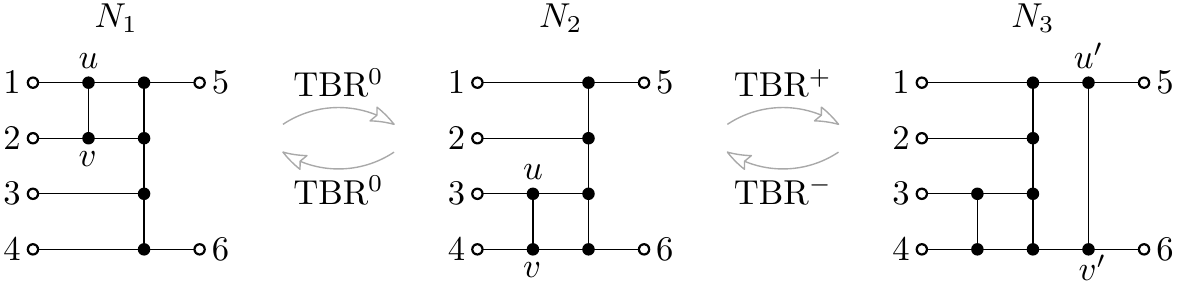}
  \caption{Illustration of the TBR operation. 
  The network $N_2$ can be obtained from $N_1$ by a \TBRZ that moves the edge $\set{u, v}$ and the network $N_3$ can be obtained from $N_2$ by a \TBRP that adds the edge $\set{u', v'}$. Each operation has its corresponding \TBRZ and \TBRM operation, respectively, that reverses the rearrangement.} 
  \label{fig:unets:TBR}
\end{figure}

\paragraph{TBR.}
A \TBR operation{\footnotemark} is the rearrangement operation that transforms a network $N\in\unets$ into another network $N' \in \unets$ in one of the following four ways: 
\begin{itemize}[leftmargin=*,label=(TBR$^-$)]
    \item[(\TBRZ)]  Remove an internal edge $e$ of $N$, subdivide an edge of the resulting graph with a new vertex $u$, subdivide an edge of the resulting graph with a new vertex $v$, and add the edge $\set{u, v}$;
   \item[ ]      or, prune an external edge $e = \set{u, v}$ of $N$ that is incident to leaf $v$ at $u$, regraft $\set{\cdot, v}$ to an edge of the resulting graph.
    \item[(\TBRP)] Subdivide an edge of $N$ with a new vertex $u$, subdivide an edge of the resulting graph with a new vertex $v$, and add the edge $e = \set{u, v}$.
    \item[(\TBRM)] Remove an edge $e$ of $N$.
\end{itemize}
\footnotetext{The TBR operation is known on unrooted phylogenetic trees as \emph{tree bisection and reconnection}.
Since in general networks are not trees and a TBR on a network does not necessarily bisect it, we use TBR now as a word on its own. For the reader who would however like to have an expansion of TBR we suggest "total branch relocation". We welcome other suggestions.}
Note that a \TBRZ can also be seen as the operation that prunes the edge $e = \set{u, v}$ at both $u$ and $v$ and then regrafts both ends. Hence, we say that a \TBRZ \emph{moves} the edge $e$. Furthermore, we say that a \TBRP \emph{adds} the edge $e$ and that a \TBRM \emph{removes} the edge $e$. These operations are illustrated in \cref{fig:unets:TBR}.  Note that a \TBRZ has an inverse \TBRZ and that a \TBRP has an inverse \TBRM, and that furthermore a \TBRP increases the reticulation number by one and a \TBRM decreases it by one.

Since a \TBR operation has to yield a phylogenetic network, there are some restrictions on the edges that can be moved or removed. Firstly, if removing an edge by a \TBRZ yields a disconnected graph, then in order to obtain a phylogenetic network an edge has to be added between the two connected components. Similarly, a \TBRM cannot remove a cut-edge. Secondly, the suppression of a vertex when removing an edge may not yield a loop $\set{u, u}$. Thirdly, removing or moving an edge cannot create a cut-edge that does not separate two leaves. Otherwise the network would not be proper.

The \TBRZ operation equals the well known TBR (tree bisection and reconnection) operation on unrooted phylogenetic trees~\cite{AS01}. The TBR operation on trees has recently been generalised to \TBRZ on improper unrooted phylogenetic networks by Francis \etal~\cite{FHMW17}.

\paragraph{PR.}
A \PR (\emph{prune and regraft}) operation is the rearrangement operation that transforms a network $N \in \unets$ into another network $N' \in \unets$ with a \PRP $=$ \TBRP, a \PRM $=$ \TBRM, or a \PRZ that prunes and regrafts an edge $e$ only at one endpoint, instead of at both like a \TBRZ. Like for TBR, we the say that the PR$^{0/+/-}$ \emph{moves/adds/removes} the edge $e$ in $N$. The PR operation is a generalisation of the well-known SPR (subtree prune and regraft) operation on unrooted phylogenetic trees~\cite{AS01}. Like for TBR, the generalisation of SPR to \PRZ for networks has been introduced by Francis \etal~\cite{FHMW17}.

\begin{figure}[htb]
  \centering
  \includegraphics{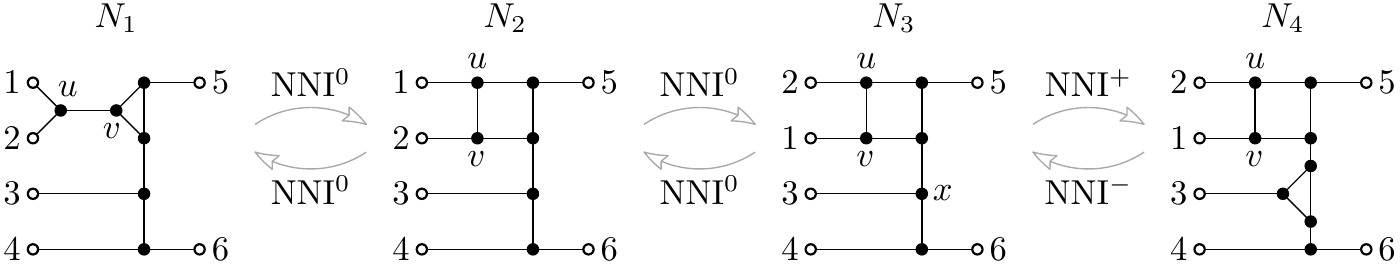}
  \caption{Illustration of the NNI operation. 
  The network $N_2$ (resp. $N_3$) can be obtained from $N_1$ (resp. $N_2$) by an \NNIZ with the axis $\{u, v\}$; alternatively, $N_2$ can be obtained from $N_1$ using the \NNIZ of $\{1,u\}$ to the triangle, and $N_3$ from $N_2$ by moving $\{1,u\}$ to the bottom edge of the square. The labels are inherited naturally following the first interpretation of the \NNIZ moves.
  The network $N_4$ can be obtained from $N_3$ by an \NNIP that extends $x$ into a triangle.
  Each operation has its corresponding \NNIZ and \NNIM operation, respectively, that reverses the transformation.} 
  \label{fig:unets:NNI}
\end{figure}

\paragraph{NNI.}
An \NNI (\emph{nearest neighbour interchange}) operation is a rearrangement operation that transforms a network $N\in \unets$ into another network $N' \in \unets$ in one of the following three ways: 
\begin{itemize}[leftmargin=*,label=(NNI$^-$)]
    \item[(\NNIZ)] Let $e= \{u, v\}$ be an internal edge of $N$. Prune an edge $f$ ($f \neq e$) at $u$, and regraft it to an edge $f'$ ($f' \neq e$) that is incident to $v$.
    \item[(\NNIP)] Subdivide two adjacent edges with new vertices $u'$ and $v'$, respectively, and add the edge $\{u', v'\}$.
    \item[(\NNIM)] If $N$ contains a triangle, remove an edge of the triangle.
\end{itemize}
These operations are illustrated in \cref{fig:unets:NNI}. We say that an \NNIZ \emph{moves} the edge $f$. Alternatively, we call the edge $e$ of an \NNIZ the \emph{axis} of the operation, as the operation can also be defined as pruning $f$ at $u$, and $f''\neq f'$ at $v$, and regrafting $f$ at $v$ and $f''$ at $u$. 
The NNI operation has been introduced on trees by Robinson~\cite{Rob71} and generalised to networks by Huber \etal~\cite{HLMW16,HMW16}.

\paragraph{Sequences and distances.}
Let $N, N' \in \unets$ be two networks.
A \emph{\TBR-sequence} from $N$ to $N'$ is a sequence
    $$\sigma = (N = N_0, N_1, N_2, \ldots, N_k = N') $$
of phylogenetic networks such that $N_i$ can be obtained from $N_{i-1}$ by a single TBR for each $i \in \set{1, 2, ..., k}$. The \emph{length} of $\sigma$ is $k$.
The \emph{\TBR-distance} $\dTBR(N, N')$ between $N$ and $N'$ is the length of a shortest TBR-sequence from $N$ to $N'$, or infinite if no such sequence exists.

Let $\class$ be a class of phylogenetic networks. 
The TBR-distance on $\class$ is defined like on $\unets$ but with the restriction that every network in a shortest TBR-sequence has to be in $\class$.
The class $\class$ is \emph{connected} under TBR if, for all pairs $N, N' \in \class$, there exists a TBR-sequence $\sigma$ from $N$ to $N'$ such that each network in $\sigma$ is in $\class$.
Hence, for the \TBR-distance to be a metric on $\class$, the class has to be connected under \TBR and the \TBR operation has to be reversible. We already noted above that the latter holds for TBR (and NNI and PR).
For a connected class $\class$, the \emph{diameter} is the maximum distance between two of its networks under its metric.
The definition for NNI and PR are analogous.

Let $\class'$ be a subclass of $\class$.
Then $\class'$ is an \emph{isometric subgraph} of a $\class$ under, say, \TBR if for every $N, N' \in \class'$ the \TBR-distance of $N$ and $N'$ in $\class'$ equals the \TBR-distance of $N$ and $N'$ in $\class$.

\section{Relations of rearrangement operations} \label{sec:relations} 
On trees, it is well known that every \NNI is also an \SPR, which, in turn, is also a \TBR.
We observe that the same holds for the generalisations of these operations as defined above.

\begin{observation} \label{clm:NNIisPRisTBR}
Let $N \in \unets$. Then, on $N$, every \NNI is a \PR and every \PR is a \TBR.
\end{observation}

For the reverse direction, we first show that every \TBR can be mimicked by at most two \PR like in $\utrees$. Then we show how to substitute a \PR with an \NNI-sequence.

\begin{lemma} \label{clm:unets:TBRisTwoPR}
Let $N, N' \in \unets$ such that $\dTBR(N, N') = 1$.
Then $1 \leq \dPR(N, N') \leq 2$, where a \TBRZ may be replaced by two \PRZ. 
  \begin{proof}
    If $N'$ can be obtained from $N$ by a \TBRP or \TBRM, then by the definition of \PRP and \PRM it follows that $\dPR(N, N') = 1$. If $N'$ can be obtained from $N$ by a \TBRZ that is also a \PRZ, the statement follows.
    Assume therefore that $N'$ can be obtained from $N$ by a \TBRZ that moves the edge $e = \set{u, v}$ of $N$ to $e' = \set{x, y}$ of $N'$. Let $G$ be the graph obtained from $N$ by removing $e$, or equivalently the graph obtained from $N'$ by removing $e'$. If $e$ is a cut-edge, then so is $e'$, and without loss of generality $u$ and $x$ as well as $v$ and $y$ subdivide an edge in the same connected components of $G$. Furthermore, if $u$ subdivides an edge of a pendant blob in $G$, then so does $x$. Otherwise $N'$ would not be proper. Therefore, the \PRZ that prunes $e$ at $u$ and regrafts it to obtain $x$ yields a phylogenetic network $N''$. The choices of $u$ and $x$ ensure that $N''$ is connected and proper. There is then a \PRZ from $N''$ to $N'$ that prunes $\set{x, v}$ at $v$ and regrafts it at $y$ to obtain $N'$. Hence, $\dPR(N, N') \leq 2$.
  \end{proof}
\end{lemma}

\begin{corollary}
Let $N, N' \in \unets$.
Then $\dTBR(N, N') \leq \dPR(N, N') \leq 2 \dTBR(N, N')$.
\end{corollary}

\begin{lemma} \label{clm:PRZtoNNIZ}
Let $N, N' \in \unetsr$ such that there is a \PRZ that transforms $N$ into $N'$. Let $e$ be the edge of $N$ pruned by this \PRZ.\\
Then there exists an \NNIZ-sequence from $N$ to $N'$ that only moves $e$ and whose length is in $\Oh(n +r)$.
Moreover, if neither $N$ nor $N'$ contains parallel edges, then neither does any intermediate networks in the \NNI-sequence.
\begin{proof}
Assume that $N$ can be transformed into $N'$ by pruning the edge $e = \{u, v\}$ at $u$ and regrafting it to $f = \{x, y\}$.
Note that there is then a (shortest) path $P = (u = v_0, v_1, v_2, \ldots, v_k = x)$ from $u$ to $x$ in $N \setminus \{e\}$, since otherwise $N'$ would be disconnected. 
Without loss of generality, assume that $P$ does not contain $y$. Furthermore, assume for now that $P$ does not contain $v$.
The idea is now to move $e$ along $P$ to $f$ with \NNIZ. In particular, we show how to construct a sequence $\sigma = (N = N_0, N_1, \ldots, N_{k} = N')$  such that either $N_{i+1}$ can be obtained from $N_{i}$ by an \NNIZ or $N_{i+1} = N_{i}$, and such that $N_i$ contains the edge $e_i = \{v_i, v\}$.
This process is illustrated in \cref{fig:unets:PRZtoNNIZ}.
Assume we have constructed the sequence up to $N_i$. 
Let $g = \{v_{i+1}, w\}$ with $w \neq v$ be the edge incident to $v_{i+1}$ that is not on $P$.
Obtain $N_{i+1}$ from $N_i$ by swapping $e_i$ and $g$ with an \NNIZ on the axis $\{v_{i}, v_{i+1}\}$.
Note that this preserves the path $P$ and that $N_{i+1}$ may only contain a parallel edge if $N$ or $N'$ contains parallel edges. As a result, we get $N_k = N'$.

\begin{figure}[htb]
  \centering
  \includegraphics{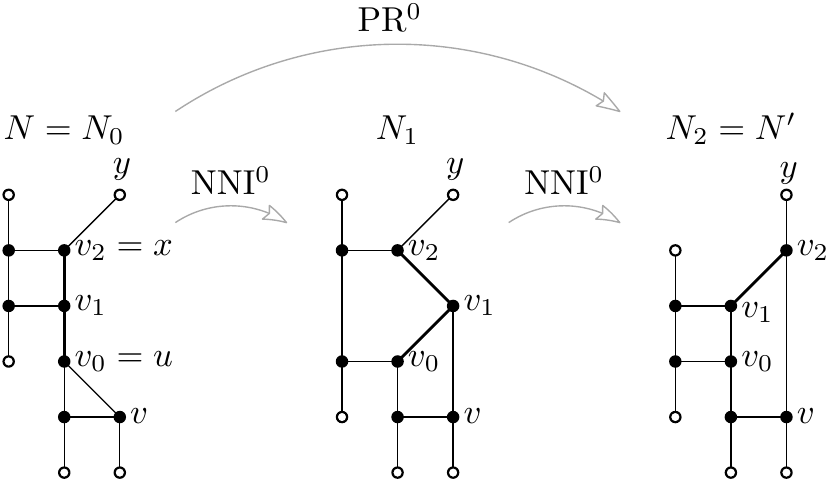}
  \caption{How to mimic the \PRZ that prunes the edge $\{u, v\}$ at $u$ and regrafts to $\{x, y\}$ with \NNIZ operations that move $u$ of $\{u, v\}$ along the path $P = (u = v_0, v_1, v_2 = x)$ (for the proof of \cref{clm:PRZtoNNIZ}). Labels follow the definition of \NNIZ along an axis.} 
  \label{fig:unets:PRZtoNNIZ}
\end{figure}

It remains to show that every network in $\sigma$ is proper. Assume otherwise and let $N_{i+1}$ be the first improper network in $\sigma$. Then $N_{i+1}$ contains a cut-edge $e_c$ that separates a blob $B$ from all leaves. We claim that $e_c$ is part of $P$. Indeed, the pruning of the \NNIZ from $N_i$ to $N_{i+1}$ has to create $B$ and the regrafting cannot be to $B$, so it has to pass along $e_c$ (\cref{fig:unets:PRZtoNNIZ:properness}). However, as $P$ is a path, the moving edge cannot pass $e_c$ again, so all networks $N_j$ for $j > i$ including $N'$ are improper; a contradiction. 
Hence, all intermediate networks $N_i$ are proper and thus $\sigma$ is an \NNIZ-sequence from $N$ to $N'$.

\begin{figure}[htb]
  \centering
  \includegraphics{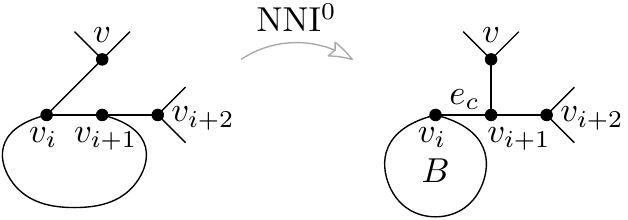}
  \caption{How an \NNIZ in the proof of \cref{clm:PRZtoNNIZ} may result an improper network where $e_c$ separates a blob $B$ from all leaves. The moving edge $\{v,v_i\}$ of $N_i$ becomes the moving edge $\{v,v_{i+1}\}$ of $N_{i+1}$. Labels follow the definition of \NNIZ along an axis.} 
  \label{fig:unets:PRZtoNNIZ:properness}
\end{figure}

Next, assume that $P$ contains $v_i = v$. Then first apply the process above to move $v$ of $\{u, v\}$ along $P' = (v = v_i, v_{i+1}, \ldots, v_k)$ to $v_k$. 
In the resulting network, apply the process above to move $u$ of $\{u, v\} = \{u, v_k\}$ along $P'' = (u = v_0, v_1, \ldots, v_i)$ to $v_i$. 
The process again avoids the creation of a network $N_j$ with parallel edges, if neither $N$ nor $N'$ contains parallel edges. Furthermore, from \cref{fig:unets:PRZtoNNIZ:properness} we get that if $\sigma$ would contain improper network then $u$ would be contained in the blob $B$. However, then $\set{u, v}$ and $e_c$ would be edges from $B$ to the rest of the network; again a contradiction.

Lastly, note that the length of $P$ is in $\Oh(n + r)$ since $N$ contains only $2n + 3r - 1$ edges. Hence, the length of $\sigma$ is also in $\Oh(n +r)$.
\end{proof}
\end{lemma}

\begin{lemma}\label{clm:PRMtoNNIM}
Let $n \geq 3$. Let $N, N' \in \unets$ such that there is a \PRM that transforms $N$ into $N'$. Let $e$ be the edge of $N$ removed by this \PRM. Let $N$ have reticulation number $r$.\\
Then, there is an \NNIZ-sequence followed by one \NNIM that transforms $N$ and $N'$ by only moving and removing $e$ and whose length is in $\Oh(n + r)$.
Moreover, if neither $N$ nor $N'$ contains parallel edges, then neither do the intermediate networks in the \NNI-sequence.
\begin{proof}
Assume the \PRM removes $e = \set{u, v}$ from $N$ to obtain $N'$. If $e$ is part of a triangle, the \PRM move is an \NNIM move. 
If $e$ is a parallel edge, then move either $u$ or $v$ with an \NNIZ to obtain a network with a triangle that contains $e$. Then the previous case applies.
So assume otherwise, namely that $e$ is not part of a triangle or a pair of parallel edges.
Then move $u$ with an \NNIZ-sequence closer to $v$ to form a triangle as follows.

Because removing $e$ in $N$ yields the proper network $N'$, it follows that $N \setminus \set{e}$ contains a shortest path $P$ from $u$ to $v$. 
Since $e$ is not part of a triangle, this path must contain at least two nodes other than  $u$ and $v$. Let $\set{x, y}$ and $\set{y, v}$ be the last two edges on $P$.
Consider the \PRZ that prunes $\{u, v\}$ at $u$ and regrafts it to $\set{x, y}$. Note that this creates a triangle on the vertices $y$, $u$ and $v$.
By \cref{clm:PRZtoNNIZ} we can replace this \PRZ with an \NNIZ-sequence. Lastly, we can remove $\{u, v\}$ with an \NNIM to obtain $N'$. The bound on the length of the \NNI-sequence as well as the second statement follow from \cref{clm:PRZtoNNIZ}.
\end{proof}
\end{lemma}

To conclude this section, we note that all previous results combined show that we can replace a \TBR-sequence with a \PR-sequence, which we can further replace with an \NNI-sequence. For several connectedness results in \cref{sec:connectedness} this allows us to focus on \TBR and then derive results for \NNI and \PR.

\section{Shortest paths} \label{sec:paths} 
In this section, we focus on bounds on the distance between two specified networks. We restrict to the \TBR-distance in $\unets$ and in $\unetsr$, and study the structure of shortest sequences of moves. We make several observations about these sequences in general, and some about shortest sequences between two networks that have certain structure in common, e.g., common displayed networks. Hence, we get bounds on the \TBR-distance between two networks, and we uncover properties of the spaces of phylogenetic networks which allow for reductions of the search space. For example, if $N$ and $N'$ have reticulation number $r$, no shortest path from $N$ to $N'$ contains a network with a reticulation number less than $r$. The proof of this statement relies on the following observation about the order in which \TBRZ and \TBRP operations can occur in a shortest path.

\begin{observation} \label{clm:unets:TBR:PMtoZ}
Let $N, N' \in \unetsr$ such that there exists a \TBR-sequence
$\sigma_0 = (N, M, N')$ that uses a \TBRP and a \TBRM. Then there is a \TBRZ that transforms $N$ into $N'$.
\end{observation}
Rephrasing \cref{clm:unets:TBR:PMtoZ}, a \TBRP followed by a \TBRM, or vice versa, can be replaced by a \TBRZ. This case can thus not occur in a shortest \TBR-sequence.
Next, we look at a \TBRZ followed by a \TBRP.

\begin{lemma} \label{clm:unets:TBR:ZPtoPZ}
Let $N, N' \in \unets$ with reticulation number $r$ and $r+1$ such that there exists a shortest \TBR-sequence $\sigma_0 = (N, M, N')$ that starts with a \TBRZ.\\
Then there is a \TBR-sequence $\sigma_{+} = (N, M', N')$ that starts with a \TBRP.
  \begin{proof}
  Note that the \TBRZ from $N$ to $M$ of $\sigma_{0}$ can be replaced with a sequence consisting of a \TBRP followed by a \TBRM. This \TBRM and the \TBRP from $M$ to $N'$ can now be combined to a \TBRZ, which gives us a sequence $\sigma_{+}$.
  \end{proof}
\end{lemma}

Let $N, N' \in \unetsr$ and consider a shortest \TBR-sequences from $N$ to $N'$ that contains \TBRP and \TBRM operations. If the reverse statement of \cref{clm:unets:TBR:ZPtoPZ} would also hold, then we could shuffle the sequence such that consecutive \TBRP and \TBRM can be replaced with a \TBRZ. This would imply that $\unetsr$ is an isometric subgraph of $\unets$ under \TBR. However, we now show that the reverse statement of \cref{clm:unets:TBR:ZPtoPZ} does not hold in general, and, hence, adjacent operations of different types in a shortest \TBR-sequence cannot always be swapped.

\begin{figure}[htb]
\centering
  \includegraphics{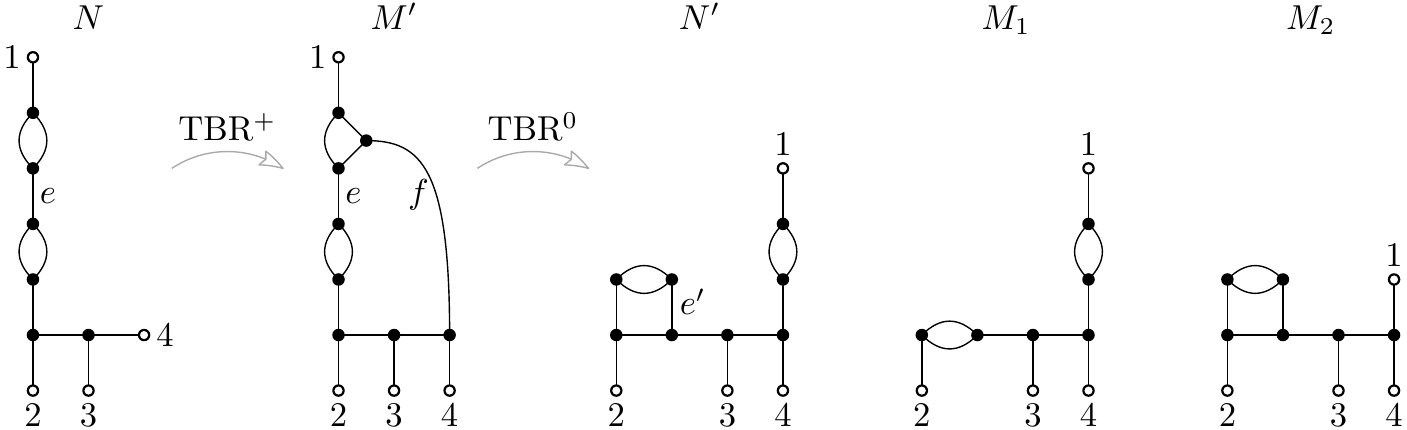}
  \caption{Two networks $N, N' \in \unets$ with TBR-distance two such that there exist a shortest TBR-sequence from $N$ to $N'$ starting with a \TBRP move (to $M'$). However, there is no shortest TBR-sequence starting with a \TBRZ, since the networks $M_1$ and $M_2$, which are the \TBRM neighbours of $N'$, have \TBRZ-distance at least two to $N$.} 
  \label{fig:unets:PZbutnoZP}
\end{figure}

\begin{lemma} \label{clm:unets:TBR:PZtoZP:notInNets}
Let $n \geq 4$ and $r \geq 2$.
Let $N, N' \in \unets$ with reticulation number $r$ and $r+1$ such that there exists a shortest \TBR-sequence
$\sigma_+ = (N, M', N')$ that starts with a \TBRP.\\
Then it is not guaranteed that there is a \TBR-sequence $\sigma_{0} =(N, M, N')$ that starts with a \TBRZ.
\begin{proof}
We claim that the networks $N$ and $N'$ in \cref{fig:unets:PZbutnoZP} are a pair of networks for which no \TBR-sequence $\sigma_{0} =(N, M, N')$ exists that starts with a \TBRZ. The two networks $M_1$ and $M_2$ in \cref{fig:unets:PZbutnoZP} are the only two \TBRM neighbours of $N'$. However, it is easy to check that the \TBRZ-distance of $N$ and $M_i$, $i \in \set{1, 2}$, is at least two. Hence, a shortest \TBR sequence from $N$ to $N'$ that starts with a \TBRZ has length three and so $\sigma_{0}$ cannot exist. Note that we can add an edge to each of the pair of parallel edges to obtain an example without parallel edges. Moreover, the example can be extended to higher $n$ and $r$ by adding extra leaves between leaf 3 and 4, and replacing a pair of parallel edges by a chain of parallel edges in each network.
\end{proof}
\end{lemma}

Note that the \TBRZ used in \cref{fig:unets:PZbutnoZP} to prove \cref{clm:unets:TBR:PZtoZP:notInNets} is a \PRZ. Hence, the statement of \cref{clm:unets:TBR:PZtoZP:notInNets} also holds for \PR. On the positive side, if one of the two networks is a tree, then we can swap the \TBRP with the \TBRZ.

\begin{figure}[htb]
  \centering
  \includegraphics{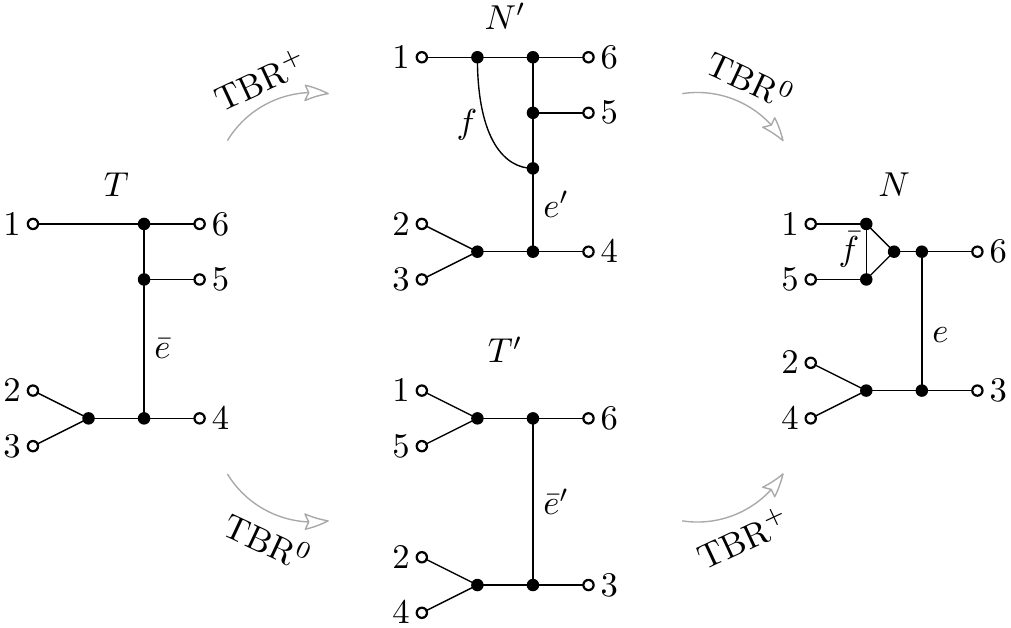}
  \caption{There is a \TBR-sequence from $T$ to $N$ that first adds $f$ with a \TBRP and then moves $e'$ to $e$ with a \TBRZ. From this, a \TBR-sequence can be derived that moves $\bar e$ to $\bar e'$ with a \TBRZ and then adds $\bar f$ with a \TBRP.} 
  \label{fig:unets:PZtoZP:trees}
\end{figure}

\begin{lemma}\label{clm:unets:TBR:PZtoZP:trees}
Let $T \in \utrees$ and $N \in \unets$ with reticulation number one such that there exists a shortest \TBR-sequence $\sigma_{+} = (T, N', N)$ that starts with a \TBRP.\\
Then there is a \TBR-sequence $\sigma_{0} =(T, T', N)$ that starts with a \TBRZ.
\begin{proof}
We show how to obtain $\sigma_{0}$ from $\sigma_{+}$. Suppose that $N'$ is obtained from $T$ by adding the edge $f$ and that $N$ is obtained from $N'$ by removing $e'$ and adding $e$. Note that $f$ is an edge of the cycle $C$ in $N'$. Furthermore, $e'$ and $f$ are distinct. Indeed, otherwise there would be a shorter \TBR-sequence from $T$ to $N$ that simply adds $e$ to $T$.

Assume for now that $e'$ is an edge of $C$ in $N'$. Then, $e'$ can be removed with a \TBRM from $N'$ to obtain a tree $T'$. Hence, the \TBRP from $T$ to $N'$ and the \TBRM from $N'$ to $T'$ can be merged into a \TBRZ from $T$ to $T'$. Furthermore, the edge $e$ can then be added to $T'$ with a \TBRP to obtain $N$. This yields the sequence $\sigma_{0}$.

Next, assume that $e'$ is not an edge of $C$ in $N'$. Then, $e'$ is a cut-edge in $N'$ and $e$ is a cut-edge in $N$. Let $\bar e$ be the edge of $T$ that equals $e'$, if it exists, or the edge that gets subdivided by $f$ into $e'$ and another edge. Let $\bar f$ be the edge of $N$ defined as follows: it is equal to $f$ itself if $f$ is not touched by the \TBRZ move from $N'$ to $N$; it is the extension of $f$ if one of its endpoints is suppressed by this move; it is one of the two edges obtained by subdividing $f$. Now let $T'$ be a tree obtained by removing $\bar f$ from $N$. Then, there is a \TBRZ from $T$ to $T'$ that moves $\bar e$ to $\bar e'$ and furthermore a \TBRP that adds $\bar f$ to $T'$ and yields $N$. We obtain again $\sigma_{0}$. An example is given in \cref{fig:unets:PZtoZP:trees}.
\end{proof}
\end{lemma}

Next, we look at shortest paths between a tree and a network. First, we show that if a network displays a tree, then there is a simple \TBRM-sequence from the network to the tree. Recall that $D(N)$ is the set of trees in $\utrees$ displayed by $N \in \unets$. This result is the unrooted analogous to Lemma 7.4 by Bordewich \etal~\cite{BLS17} on rooted phylogenetic networks.

\begin{lemma} \label{clm:unets:TBR:pathDown}
Let $N \in \unetsr$ and $T \in \utrees$. \\
Then $T \in D(N)$ if and only if $\dTBR(T, N) = r$, that is, iff there exists a \TBRM-sequence of length $r$ from $N$ to $T$.
\begin{proof}
Note that $\dTBR(T, N) \geq r$, since a \TBR can reduce the reticulation number by at most one. Furthermore, if we apply a sequence of $r$ \TBRM moves on $N$, we arrive at a tree that is displayed by $N$. Hence, if $T \not\in D(N)$, then $\dTBR(T,N) > r$.

We now use induction on $r$ to show that $\dTBR(T, N) \leq r$ if $T\in D(N)$. If $r = 0$, then $T = N$ and the inequality holds. Now suppose that $r > 0$ and that the statement holds whenever a network with a reticulation number less than $r$ displays $T$. Fix an embedding of $T$ into $N$ and colour all edges of $N$ not covered by this embedding green. Note that removing a green edge with a \TBRM might result in an improper network or a loop. Therefore, we have to show that there is always at least one edge that can be removed such that the resulting graph is a phylogenetic network. For this, consider the subgraph $H$ of $N$ induced by the green edges. If $H$ contains a component consisting of a single green edge $e$, then removing $e$ from $N$ with a \TBRM yields a network $N'$. If $H$  contains a tree component $S$, then it is easy to see that removing an external edge of $S$ from $N$ with a \TBRM yields a network $N'$. Otherwise, as $N$ is proper, a component $S$ displays a tree $T_S$ whose external edges cover exactly the external edges of $S$.
We can then apply the same case distinction to the edges of $S$ not covered by $T_S$ and either directly find an edge to remove or find further trees that cover the smaller remaining components. Since $S$ is finite, we eventually find an edge to remove. The induction hypothesis then applies to $N'$. This concludes the proof.
\end{proof} 
\end{lemma}

Note that the proof of \cref{clm:unets:TBR:pathDown} also works if $T$ is a network displayed by $N$. Hence, we get the following corollary.

\begin{corollary} \label{clm:unets:TBR:pathDown:nets}
Let $N \in \unetsr$ and let $N' \in \unetsxx[n,r']$ such that $N'$ is displayed by $N$.\\
Then $\dTBR(N', N) = r - r'$, that is, there exists a \TBRM-sequence of length $r-r'$ from $N$ to $N'$.
\end{corollary}

\Cref{clm:unets:TBR:pathDown} and \cref{clm:unets:TBR:pathDown:nets} now allow us to construct \TBR-sequences between networks that go down tiers and then come up again. In fact, for rooted networks this can sometimes be necessary as Klawitter and Linz have shown~\cite[Lemma 13]{KL19}. However, we now show that this is never necessary for \TBR on unrooted networks.

\begin{lemma} \label{clm:unets:TBR:noNeedToGoDown}
Let $N, N' \in \unets$.\\ 
Then in no shortest \TBR-sequence from $N$ to $N'$ does a \TBRM precede a \TBRP.
  \begin{proof}
    Consider a minimal counterexample with $N, N' \in \unets$ such that there exists a shortest \TBR-sequence $\sigma$ from $N$ to $N'$ that uses exactly one \TBRM and \TBRP and that starts with this \TBRM. If $\sigma$ uses \TBRZ operations between the \TBRM and the \TBRP, then, by \cref{clm:unets:TBR:ZPtoPZ}, we can swap the \TBRP forward until it directly follows the \TBRM. However, then we can obtain a \TBR-sequence shorter than $\sigma$ by combining the \TBRM and the \TBRP into a \TBRZ by \cref{clm:unets:TBR:PMtoZ}; a contradiction.
  \end{proof}
\end{lemma}

Combining \cref{clm:unets:TBR:pathDown,clm:unets:TBR:pathDown:nets,clm:unets:TBR:ZPtoPZ}, we easily derive the following two corollaries about short sequences that do not go down tiers before going back up again.

\begin{corollary} \label{clm:unets:TBR:distanceViaDisplayedTrees}
Let $N, N' \in \unets$ with reticulation number $r$  and $r'$, with $r\geq r'$. Then
    $$\dTBR(N, N') \leq \min\set{\dTBR(T, T') \colon T \in D(N), T' \in D(N')} + r \text{.}$$
\end{corollary}

\begin{corollary} \label{clm:unets:TBR:distanceSharedDisplayedTrees}
Let $N, N' \in \unets$ with reticulation number $r$  and $r'$, and $r\geq r'$. Let $T \in \utrees$ such that $T \in D(N), D(N')$. Then
    $$\dTBR(N, N') \leq r \text{.}$$
\end{corollary}

Both \cref{clm:unets:TBR:distanceViaDisplayedTrees,clm:unets:TBR:distanceSharedDisplayedTrees} can easily be proven by first finding a sequence that goes down to tier 0 and back up to tier $r$, and then combining the $r'$ \TBRM with $r'$ \TBRP into $r'$ \TBRZ using \cref{clm:unets:TBR:ZPtoPZ}.

The following lemma is the unrooted analogue to Proposition 7.7 by Bordewich\linebreak \etal~\cite{BLS17}. We closely follow their proof.
\begin{lemma} \label{clm:unets:TBR:existingCloseDisplayedTree}
Let $N, N' \in \unets$ such that $\dTBR(N, N') = k$. Let $T \in D(N)$. \\
Then there exists a $T' \in D(N)$ such that
    $$\dTBR(T, T') \leq k \text{.}$$
  \begin{proof}
  The proof is by induction on $k$. If $k = 0$, then the statement trivially holds. Suppose that $k = 1$. If $T \in D(N')$, then set $T' = T$, and we have $\dTBR(T, T') = 0 \leq 1$. So assume otherwise, namely that $T \not \in D(N')$. Note that that if $N'$ has been obtained from $N$ by a \TBRP, then $N'$ displays $T$. Therefore, distinguish whether $N'$ has been obtained from $N$ by a \TBRZ or \TBRM $\sigma$.
  
  Suppose that $N'$ has been obtained from $N$ by a \TBRZ that moves the edge $e = \set{u, v}$ of $N$. Fix an embedding $S$ of $T$ into $N$. Since $N'$ does not display $T$, the edge $e$ is covered by $S$. Let $\bar e$ be the edge of $T$ that gets mapped to the path of $S$ that covers $e$. Let $S_1$ and $S_2$ be the subgraphs of $S \setminus \set{e}$. Note that $S_1, S_2$ have embeddings into $N$ and $N'$. Now, if in $N$ there exists a path $P$ from the embedding of $S_1$ to the embedding of $S_2$ that avoids $e$, then the graph consisting of $P$, $S_1$, and $S_2$ is a tree $T'$ displayed by $N'$. Otherwise $e$ is a cut-edge of $N$ and the \TBRZ moves $e$ to an edge $e'$ connecting the two components of $N \setminus \set{e}$. Then in $N'$ there is path $P$ from the embedding of $S_1$ to the embedding of $S_2$ in $N'$. Together they form an embedding of a tree $T'$ displayed by $N'$. In both cases $T'$ can also be obtained from $T$ by moving $\bar e$ to where $P$ attaches to $S_1$ and $S_2$. If $N'$ is obtained from $N$ by a \TBRM, then the first case has to apply.
  
  Now suppose that $k \geq 2$ and that the hypothesis holds for any two networks with \TBR-distance at most $k-1$. Let $N'' \in \unets$  such that $\dTBR(N, N'') = k-1$ and $\dTBR(N'', N') = 1$. Thus by induction there are trees $T''$ and $T'$ such that $T'' \in D(N'')$ with $\dTBR(T, T'') \leq k-1$ and $T' \in D(N')$ with $\dTBR(T'', T') \leq 1$. It follows that $\dTBR(T, T') \leq k$, thereby completing the proof of the lemma.
  \end{proof}
\end{lemma}

By setting one of the two networks in the previous lemma to be a phylogenetic tree and noting that the roles of $N$ and $N'$ are interchangeable, the next two corollaries are immediate consequences of \cref{clm:unets:TBR:pathDown,clm:unets:TBR:existingCloseDisplayedTree}.

\begin{corollary} \label{clm:unets:TBR:distanceToDisplayedTree}
Let $T \in \utrees$, $N \in \unetsr$ such that $\dTBR(T, N) = k$. 
Then for every $T' \in D(N)$ 
    $$\dTBR(T, T') \leq k \text{.}$$
\end{corollary}

\begin{corollary} \label{clm:unets:TBR:distanceOfDisplayedTrees}
Let $N \in \unetsr$ and let $T, T' \in D(N)$. 
Then
    $$\dTBR(T, T') \leq r \text{.}$$
\end{corollary}

The following theorem is the unrooted analogous of Theorem 7 by Klawitter and Linz~\cite{KL19} and their proof can be applied straightforward by swapping SNPR and rooted networks with TBR and unrooted networks, and by using \cref{clm:unets:TBR:pathDown,clm:unets:TBR:existingCloseDisplayedTree} and \cref{clm:unets:TBR:treesIsometric}.
\begin{theorem} \label{clm:unets:TBR:distanceTreeNetwork}
Let $T \in \utrees$ and let $N \in \unetsr$. Then
    $$\dTBR(T, N) = \min\limits_{T' \in D(N)} \dTBR(T, T') + r \text{.}$$
\end{theorem}

\section{Connectedness and diameters} \label{sec:connectedness} 
Whereas in the previous section we studied the maximum distance between two given networks, here, we focus on global connectivity properties of several classes of phylogenetic networks under NNI, PR, and TBR. These results imply that these operations induce metrics on these spaces. For each connected metric space, we can ask about its diameter. Since a class of phylogenetic networks that contains networks with unbounded reticulation number naturally has an unbounded diameter, this questions is mainly of interest for the tiers of a class. First, we recall some known results from unrooted phylogenetic trees.
\begin{theorem}[Li \etal \cite{LTZ96}, Ding \etal \cite{DGH11}] \label{clm:utrees:diameter}
The space $\utrees$ is connected under 
\begin{itemize}
    \item \NNIZ with the diameter in $\Theta(n \log n)$,
    \item \PRZ with the diameter in $n - \Theta(\sqrt{n})$, and
    \item \TBRZ with the diameter in $n - \Theta(\sqrt{n})$.
\end{itemize}
\end{theorem}

\subsection{Network space}
Huber \etal~\cite[Theorem 5]{HMW16} proved that the space of phylogenetic networks that includes improper networks is connected under \NNI. We reprove this for our definition of $\unets$, but first look at the tiers of this space.

\begin{theorem} \label{clm:unets:NNIconnected:tier}
Let $n \geq 0$, $r \geq 0$, and $m = n + r$.\\
Then $\unetsr$ is connected under \NNI with the diameter in $\Theta(m \log m)$.
\begin{proof}
Let $N \in \unetsr$ and let $T \in \utrees$ be a tree displayed by $N$. 
We show that $N$ can be transformed into a sorted $r$-handcuffed caterpillar $N^*$ with $\Oh(m \log m)$ \NNI. Our process is as follows and illustrated in \cref{fig:unets:NNIdiam:process}.
\begin{description} 
    \item[Step 1.] Transform $N$ into a network $N_T$ that is tree-based on $T$.
    \item[Step 2.] Transform $N_T$ into handcuffed tree $N_H$ on the leafs 1 and 2.
    \item[Step 3.] Transform $N_H$ into a sorted handcuffed caterpillar $N^*$.
\end{description}
\begin{figure}[htb]
  \centering
  \includegraphics{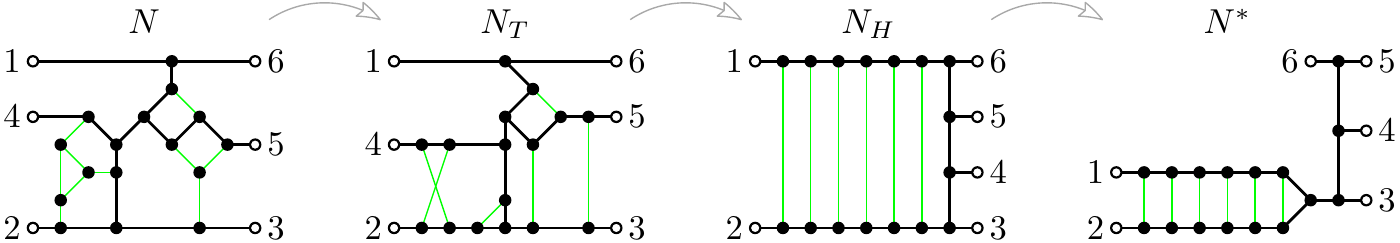}
  \caption{The process used in the proof of \cref{clm:unets:NNIconnected:tier}.
  We transform a network $N$ into a tree-based network $N_T$, then into a handcuffed tree $N_H$, and finally into a sorted handcuffed caterpillar $N^*$.}
  \label{fig:unets:NNIdiam:process}
\end{figure}

We now describe this process in detail. 
For \textbf{Step 1}, we show how to construct an \NNIZ-sequence $\sigma$ from $N$ to $N_T$, and we give a bound on the length of $\sigma$.
Let $S$ be an embedding of $T$ into $N$, that is, $S$ is a subdivision of $T$ and a subgraph of $N$. Colour all edges of $N$ used by $S$ black and all other edges green.
Note that this yields green, connected subgraphs $G_1, \ldots, G_l$ of $N$; more precisely, the $G_i$ are the connected components of the graph induced by the green edges of $N$.
Note that each $G_i$ has at least two vertices in $S$, since otherwise $N$ would not be proper.
Furthermore, if each $G_i$ consists of a single edge, then $N$ is tree-based on $T$.
Assuming otherwise, we show how to break the $G_i$ apart.

First, if there is a triangle on vertices $v_1, u, v_2$ where $v_1$ and $v_2$ are 
adjacent vertices in $S$ and $u$ is their neighbour in $G_i$, then change the embedding of $S$ (and $T$) so that it takes the path $v_1, u, v_2$ instead of $v_1, v_2$ (see \cref{fig:unets:NNIdiam:tbased}a).
Otherwise, there is an edge $\set{v, u}$ where $v$ is in $S$ and the other vertices adjacent to $u$ are not adjacent to $v$. Let $\set{u, w_1}$ and $\set{u, w_2}$ be the other edges incident to $u$. Apply an \NNIZ to move $\set{u, w_1}$ to $S$ as in \cref{fig:unets:NNIdiam:tbased}b.
Note that each such \NNIZ decreases the number of vertices in green subgraphs and increases the number of vertices in $S$. Furthermore, the resulting networks is clearly proper.
Therefore, repeat these cases until all $G_i$ consist of single edges.
Let the resulting graph be $N_T$.
Since there are at most $2(r-1)$ vertices in all green subgraphs that are not in $S$, the number of required \NNIZ for Step 1 is at most
\begin{equation} \label{eq:unets:NNI:diam1}
    2(r-1)\text{.}
\end{equation}

\begin{figure}[htb]
  \centering
  \includegraphics{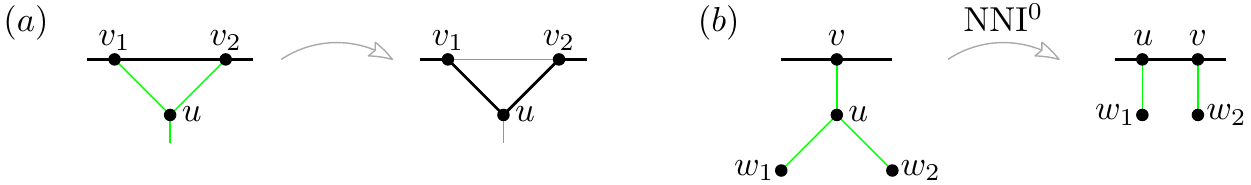}
  \caption{Transformation and \NNIZ used in Step 1 to obtain a tree-based network $N_T$.}
  \label{fig:unets:NNIdiam:tbased}
\end{figure}

In \textbf{Step 2} we transform $N_T$ into a handcuffed tree $N_H$ on the leaves 1 and 2.
Let $M = \set{\set{u_1, v_1}, \set{u_2, v_2}, \ldots, \set{u_r, v_r}}$ be the set of green edges in $N_T$, that is, the edges that are not in the embedding $S$ of $T$ into $N_T$.  
Without loss of generality, assume that for $i \in \set{1, \ldots, r}$ the distance between $u_i$ and leaf $1$ in $S$ is at most the distance of $v_i$ to leaf $1$ in $S$.
The idea is to sweep along the edges of $S$ to move the $u_i$ towards leaf $1$ and then do the same for the $v_i$ towards leaf $2$.

For an edge $e$ of $T$, let $P_e$ be the path of $S$ corresponding to $e$.
Let $e_1$ be the edge of $T$ incident to leaf $1$.
Impose directions on the edges of $T$ towards leaf $1$. Do the same for the edges of $S$ accordingly.
This gives a partial order $\preceq$ on the edges of $T$ with $e_1$ as maximum.
Let $\prec$ be a linear extension of $\preceq$ on the edges of $T$.

Let $e = (x, y)$ be the minimum of $\prec$. 
Let $P_e = (x, \ldots, y)$ be the corresponding path in $S$.
From $x$ to $y$ along $P_e$, proceed as follows.
\begin{enumerate}[label=(\roman*)]
    \item If there is an edge $(u_i, v_l)$ in $P_e$, then swap $u_i$ and $v_l$ with an \NNIZ. 
    \item If there is an edge $(u_i, u_j)$ in $P_e$ then move the $u_j$ endpoint of the green edge incident to $u_j$ onto the green edge incident to $u_i$ with an \NNIZ.
    \item Otherwise, if there is an edge $(u_i, y)$ in $P_e$, then move $u_i$ beyond $y$.
\end{enumerate}
This is illustrated in \cref{fig:unets:NNIdiam:sweep}. Informally speaking, we stack $u_j$ onto $u_i$ so they can move together towards $e_1$.
Repeat this process for each edge in the order given by $\prec$.
For the last edge $e_1$, ignore case (iii).
Next ``unpack'' the stacked $u_i$'s on $e_1$. 

We now count the number of \NNIZ needed.
Firstly, each $v_l$ is swapped at most once with a $u_i$. 
Secondly, each $u_j$ is moving to and from a green edge at most once.
Furthermore, each vertex of $S$ corresponding to a vertex of $T$ is swapped at most twice. 
Hence, the total number of \NNIZ required is at most 
\begin{equation} \label{eq:unets:NNI:diam2}
    3r + 2n \text{.}
\end{equation}

\begin{figure}[htb]
  \centering
  \includegraphics{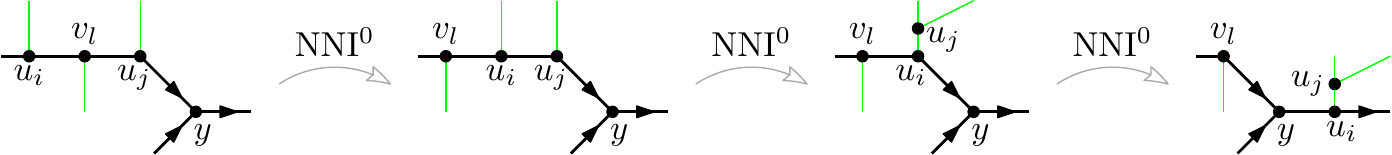}
  \caption{\NNIZ used in Step 2 
  to obtain a hand-cuffed tree $N_H$. The label of the moving endpoint follows this endpoint to its regrafting point.}
  \label{fig:unets:NNIdiam:sweep}
\end{figure}

Repeat this process for the $v_i$ towards leaf $2$. Since the $v_i$ do not have to be swapped with $u_j$, the total number of \NNIZ required for this is at most
\begin{equation} \label{eq:unets:NNI:diam3}
    2r + 2n \text{.}
\end{equation}
Note that the resulting network may not yet be a handcuffed tree as the order of the $u_i$ and $v_j$ may be different.
Hence, lastly in Step 2, to obtain $N_H$ sort the edges with the mergesort-like algorithm by Li \etal~\cite[Lemma 2]{LTZ96}. They show that the required number of \NNIZ for this is at most
\begin{equation} \label{eq:unets:NNI:diam4}
    r (1 + \log r) \text{.}
\end{equation}

For \textbf{Step 3}, consider the path $P$ in $S$ from leaf $1$ to $2$. 
If $P$ contains only one pendant subtree, then $N_H$ is handcuffed on the cherry $\set{1, 2}$. 
Otherwise, use \NNIZ to reduce it to one pendant subtree. This takes at most $n$ \NNIZ.
Next, transform the pendant subtree of $P$ into a caterpillar to obtain a handcuffed caterpillar, again with at most $n$ \NNIZ.
Lastly, sort the leaves with the algorithm from Li \etal~\cite[Lemma 2]{LTZ96} to obtain the sorted handcuffed caterpillar $N^*$.
The required number of \NNIZ to get from $N_H$ to $N^*$ is at most
\begin{equation} \label{eq:unets:NNI:diam5}
    2n + n \log n \text{.}
\end{equation}

Since we can transform any network $N\in\unetsr$ into $N^*$, it follows that $\unetsr$ is connected under $\NNI$.
Furthermore, adding \crefrange{eq:unets:NNI:diam1}{eq:unets:NNI:diam5} up and multiplying the result by two shows that the diameter of $\unetsr$ under \NNIZ is at most 
\begin{equation} \label{eq:unets:NNI:diam6}
    2(6n + 8r + n \log n + r \log r) \in \Oh((n + r) \log (n + r)) \text{.}
\end{equation}
Francis \etal~\cite[Theorem 2]{FHMW17} gave the lower bound $\Omega(m \log m)$ on the diameter of tier $r$ of the space that allows improper networks under \NNIZimproper (\NNIZ without the properness condition). 
Their proof consists of two parts: a lower bound on the total number of networks in a tier $\abs{\unetsr}$, and upper bounds on the number of networks that can be reached from one network for each fixed number of \NNIZimproper.
The diameter of $\unetsr$ is at least the smallest number of moves needed for which previously mentioned upper bound is greater than the lower bound on $\abs{\unetsr}$. 

Our version of \NNIZ is stricter than theirs as we do not allow improper networks. Hence, the number of networks that can be reached with a fixed number of \NNIZ is at most the number of networks that can be reached with the same number of \NNIZimproper. Furthermore, their lower bound on $\abs{\unetsr}$ is found by counting the number of \emph{Echidna} networks, a class of networks only containing proper networks. Combining these two observations, we see that their lower bound for the diameter of $\unetsr$ under \NNIZimproper is also a lower bound for $\unetsr$ under \NNIZ.
\end{proof}
\end{theorem}

From \cref{clm:unets:NNIconnected:tier} we get the following corollary.
\begin{corollary} \label{clm:unets:NNI:connected}
The space $\unets$ is connected under \NNI with unbounded diameter.
\end{corollary}

Since, by \cref{clm:NNIisPRisTBR}, every \NNI is also a \PR and \TBR, the statements in  \cref{clm:unets:NNIconnected:tier} and \cref{clm:unets:NNI:connected} also hold for \PR and \TBR. This observation has been made before by Francis \etal \cite{FHMW17} for tiers of the space of networks that allow improper networks. 
\begin{corollary} \label{clm:unets:TBR:connected}
The spaces $\unets$ and $\unetsr$ are connected under the \PR and \TBR operation.
\end{corollary}

We now look at the diameters of $\unetsr$ under \PR and \TBR.
\begin{theorem} \label{clm:unets:PR:tierConnected}
Let $n \geq 0$, $r \geq 0$.\\
Then the diameter of $\unetsr$ under \PRZ is in $\Theta(n + r)$ with the upper bound $n + 2r$.
\begin{proof}
The asymptotic lower bound was proven by Francis \etal~\cite[Proposition 4]{FHMW17}.
Concerning an upper bound, Janssen \etal~\cite[Theorem 4.22]{JJEvIS17} showed that the distance of two improper networks $M$ and $M'$ under \PR is at most $n + \frac{8}{3}r$, of which $\frac{2}{3}r$ \PRZ moves are used to transform $M$ and $M'$ into proper networks $N$ and $N'$. Hence, the \PR-distance of $N$ and $N'$ is at most $n + 2r$.
\end{proof}
\end{theorem}

\begin{theorem} \label{clm:unets:TBR:tierConnected}
Let $n \geq 0$, $r \geq 0$.\\
Then the diameter of $\unetsr$ under \TBR is in $\Theta(n + r)$ with the upper bound $$n - 3 - \floor{\frac{\sqrt{n - 2} - 1}{2}} + r\text{.}$$
\begin{proof}
Like for PR, the lower bound was proven by Francis \etal~\cite[Proposition 4]{FHMW17}. In \cref{clm:unets:TBR:distanceViaDisplayedTrees} we show that the \TBR-distance of two networks $N$ and $N' \in \unetsr$ that display a tree $T$ and $T' \in \utrees$, respectively, is at most $\dTBR(T, T') + r$. Since $\dTBR(T, T') \leq n - 3 - \floor{\frac{\sqrt{n - 2} - 1}{2}}$ by Theorem 1.1 of Ding \etal~\cite{DGH11} it follows that $\dTBR(N, N') \leq n - 3 - \floor{\frac{\sqrt{n - 2} - 1}{2}} + r$.
\end{proof}
\end{theorem}

\subsection{Networks displaying networks}
Bordewich~\cite[Proposition 2.9]{Bor03} and Mark \etal~\cite{MMS16} showed that the space of rooted phylogenetic trees that display a set of triplets (trees on three leaves) is connected under \NNI. 
Furthermore, Bordewich \etal~\cite{BLS17} showed that the space of rooted phylogenetic networks that display a set of rooted phylogenetic trees is connected.
We give a general result for unrooted phylogenetic networks that display a set of networks.
For this, we will use \cref{clm:unets:TBR:pathDown}, which, as we recall, guarantees that if a network $N \in \unetsr$ displays a tree $T \in \utrees$, then there is a sequence of $r$ \TBRM from $N$ to $T$. 

\begin{proposition} \label{clm:unets:displayingSet:connectivity}
Let $P = \set{P_1, ..., P_k}$ be a set of $k$ phylogenetic networks $P_i$ on $Y_i \subseteq X = \set{1, \ldots, n}$.\\
Then $\unets(P)$ is connected under \NNI, \PR, and \TBR.
\begin{proof}
Define the network $N_P \in \unets(P)$ as follows. Let $P_0 \in \utrees$ be the caterpillar where the leaves are ordered from $1$ to $n$; that is, $P_0$ contains a path $(v_2, v_3, \ldots, v_{n-1})$ such that leaf $i$ is incident to $v_i$, leaf $1$ is incident to $v_2$, and leaf $n$ is incident to $v_{n-1}$. Let $e_i$ be the edge incident to leaf $i$ in $P_0$.
Subdivide $e_i$ with $k$ vertices $u_i^1, \ldots, u_i^k$.
Now, for $P_j \in P$, $j \in \set{1, \ldots, k}$, identify leaf $i$ of $P_j$ with $u_i^j$ of $P_0$ and remove its label $i$. 
Finally, in the resulting network suppress any degree two vertex. This is necessary if one or more of the $P_j$ have fewer than $n$ leaves. 
The resulting network $N_P$ now displays all networks in $P$. An example is given in \cref{fig:unets:CanonicalDisplayingNetwork}.

\begin{figure}[htb]
\begin{center}
  \includegraphics{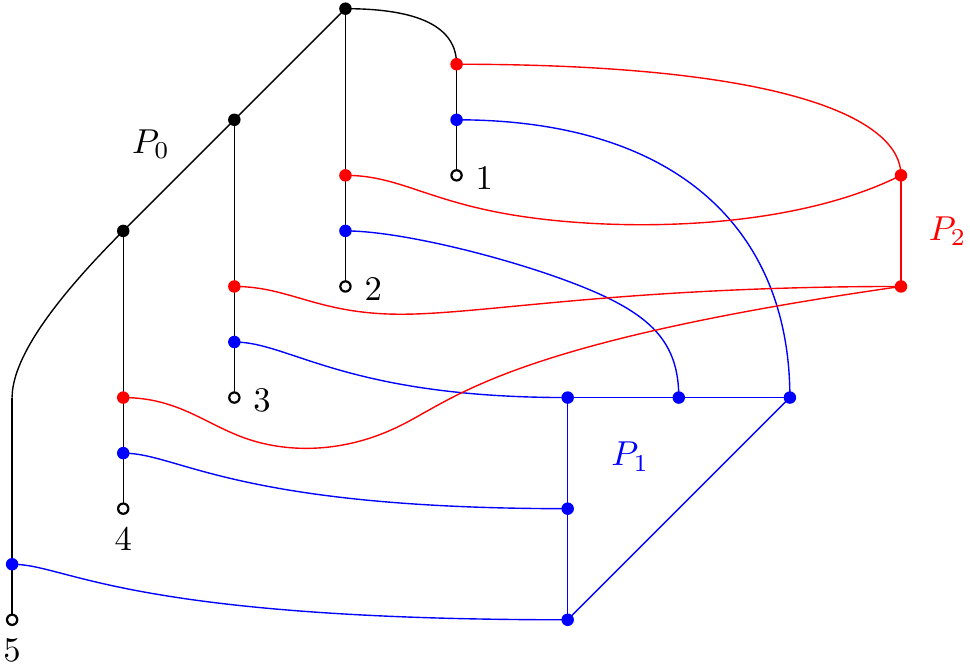}
  \caption{The canonical network $N_P \in \unetsx[5]$ that displays the set of phylogenetic networks $P = (P_1, P_2)$ with the underlying caterpillar $P_0$.} 
  \label{fig:unets:CanonicalDisplayingNetwork}
\end{center}
\end{figure}

Let $N \in \unets(P)$. Construct a \TBR-sequence from $N$ to $N_P$ by, roughly speaking, building a copy of $N_P$ attached to $N$, and then removing the original parts of $N$. First, add $P_0$ to $N$ by adding an edge $e = \set{v_1, v_2}$ from the edge incident to leaf 1 to the edge incident to leaf 2 with a \TBRP. Then add another edge from $e$ to the edge incident to leaf 3, and so on up to leaf $n$. Colour all newly added edges and the edges incident to the leaves blue, and all other edges red. Note that the blue edges now give an embedding of $P_0$ into the current network. Now, ignoring all red edges, it is straight forward to add the $P_j$, $j \in \{1, \ldots, k\}$ one after the other with \TBRP such that the resulting network displays $N_P$. For example, one could start by adding a tree displayed by $P_j$ and then adding any other edges.  The first part works similar to the construction of $P_0$ and the second part is possible by \cref{clm:unets:TBR:pathDown}. Lastly, remove all red edges with \TBRM such that every intermediate network is proper. This is again possible by \cref{clm:unets:TBR:pathDown} and yields the network $N_P$. Note that in the first two stages the red edges (plus external edges) display $P$ and in the last phase the non-red edges display $P$.

Since we only used \TBRP and \TBRM operations, the statement also holds for \PR. For \NNI, by \cref{clm:PRMtoNNIM} we can replace each of these operations that add or remove an edge $e$ by \NNI-sequences that only move and remove or add the edge $e$. Hence, the statement also holds for \NNI.
\end{proof}
\end{proposition}

For the following corollary, note that a quartet is an unrooted binary tree on four leaves and a quarnet is an unrooted binary, level-1 network on four leaves \cite{HMSW18}.

\begin{corollary}
Let $X = \set{1, ..., n}$.
Let $P$ be a set of phylogenetic trees on $X$, a set of quartets on $X$, or a set of quarnets on $X$.
Then $\unets(P)$ is connected under \NNI, \PR, and \TBR.
\end{corollary}

\subsection{Tree-based networks}
A related but more restrictive concept to displaying a tree is being tree-based. So, next, we consider the class of tree-based networks. We start with the tiers of $\utbasednets(T)$, which is the set of tree-based networks that have the tree $T$ as base tree.

\begin{theorem} \label{clm:unets:tbased:connectedness}
Let $T \in \utrees$.
Then the space $\utbasednetsr(T)$ is connected under
\begin{itemize}
    \item \TBR with the diameter being between $\ceil{\frac{r}{3}}$ and $r$,
    \item \PR with the diameter being between $\ceil{\frac{r}{2}}$ and $2r$, and
    \item \NNI with the diameter being in $\Oh(r(n + r))$.
\end{itemize}
\begin{proof}
We start with the proof for TBR. Let $N, N' \in \utbasednetsr(T)$. Consider embeddings of $T$ into $N$ and $N'$. Let $S = \set{e_1, \ldots, e_r}$ and $S' = \set{e_1', \ldots, e_r'}$ be the set of all edges not covered by this embedding of $T$ in $N$ and in $N'$. Since $N$ is tree-based, $S$ and $S'$ consist of vertex-disjoint edges. Following the embeddings of $T$ into $N$ and $N'$, it is straightforward to move each edge $e_i$ with a \TBRZ from $N$ to where $e_i'$ is in $N'$. In total, this requires at most $r$ \TBRZ. Since every intermediate network is clearly in $\utbasednetsr(T)$, this gives connectedness of $\utbasednetsr(T)$ and an upper bound of $r$ on the diameter. For the lower bound, consider a network $M$ with $r$ pairs of parallel edges and $M'$ without any. Observe that a \TBRZ can break at most three pairs of parallel edges and that only if a pair of parallel edges is removed and attached to two other pairs of parallel edge. Hence, for these particular $N$ and $N'$ we have that $\dTBR(N, N') \geq \ceil{\frac{r}{3}}$. 

The constructed \TBRZ-sequence for $N$ to $N'$ above can be converted straightforwardly into a \PRZ-sequence from $N$ to $N'$ of length at most $2r$. For the lower bound, let $M$ and $M'$ be as above and note that a \PR can break at most two pairs of parallel edges. Hence, $\dPR(M, M') \geq \ceil{\frac{r}{2}}$.

By \cref{clm:PRZtoNNIZ}, the \PR-sequence can be used to construct an \NNI-sequence from $N$ to $N'$ that only moves the edges $e_i$ along paths of the embedding of $T$. Since the \PR-sequence has length at most $2r$ and each \PR can be replaced by an \NNI sequence of length at most $\Oh(n + r)$, this gives the upper bound of $\Oh(r(n + r))$ on the diameter of $\utbasednetsr(T)$ under \NNI. 
\end{proof}
\end{theorem} 

We use \cref{clm:unets:tbased:connectedness} to prove connectedness of other spaces of tree-based networks.

\begin{theorem} \label{clm:unets:tbased:connectednessOther}
Let $T \in \utrees$.\\
Then the spaces $\utbasednets(T)$, $\utbasednetsr$, and $\utbasednets$ are each connected under \TBR, \PR, and \NNI.
Moreover, the diameter of $\utbasednetsr$ is in $\Theta(n + r)$ under \TBR and \PR and in $\Oh(n \log n + r(n + r))$ under \NNI.
\begin{proof}
Assume without loss of generality that $T$ has the cherry $\set{1, 2}$.
First, let $N$ and $N'$ be in tiers $r$ and $r'$ of $\utbasednets(T)$, respectively, such that they are $r$- and $r'$-handcuffed on the cherry $\set{1, 2}$. Then $\dNNI(N, N') = \abs{r' - r}$, as we can decrease the number of handcuffs with \NNIM. 
Since, by \cref{clm:unets:tbased:connectedness}, the tiers of $\utbasednetsr(T)$ are connected, the connectedness of $\utbasednets(T)$ follows.

Second, let $N, N' \in \utbasednetsr$ be tree-based networks on $T$ and $T'$ respectively, and with an $r$-burl on the edge incident to leaf $1$. 
Ignoring the burls, by \cref{clm:utrees:diameter}, $N$ can be transformed into $N'$ by transforming $T$ into $T'$ with $\Oh(n \log n)$ \NNIZ or with $\Oh(n)$ \PRZ or \TBRZ. With \cref{clm:unets:tbased:connectedness}, the connectedness of $\utbasednetsr$ and the upper bounds on the diameter follow.
The lower bound on the diameter under \PR and \TBR also follows from \cref{clm:utrees:diameter} and \cref{clm:unets:tbased:connectedness}, 

Lastly, the connectedness of $\utbasednets$ follows similarly from the connectedness of $\utrees$ and $\utbasednetsr$.
\end{proof}
\end{theorem}

\subsection{Level-$k$ networks}
To conclude this section, we prove the connectedness of the space of level-$k$ networks.

\begin{theorem}\label{clm:unets:lvlk:connectedness:TBR}
Let $n \geq 2$ and $k \geq 1$.\\
Then, the space $\ulvlknets$ is connected under \TBR and \PR with unbounded diameter.
\begin{proof}
Let $N \in \ulvlknets$ and $T \in \utrees$.
We show that $N$ can be transformed into the network $M \in \ulvlknets$ that can be obtained from $T$ by adding a $k$-burl to the edge incident to leaf $1$. 
First, create a $k$-burl in $N$ on the edge incident to leaf $1$. This can be done using $k$ \PRP.
Next, using  \cref{clm:unets:TBR:pathDown} remove all other blobs. This gives a network $M'$ which consists of a tree $T'$ with a $k$-burl at leaf $1$. There is a \PRZ-sequence from $T'$ to $T$, which is easily converted into a sequence from $M'$ to $M$. This proves the connectedness of $\ulvlknets$ under \PR and also \TBR.
Lastly, note that the diameter is unbounded because the number of possible reticulations in a level-$k$ network is unbounded.
\end{proof}
\end{theorem}

Note that an \NNIP cannot directly create a pair of parallel edges. 
We may instead add a triangle with an \NNIP and then use an \NNIZ to transform it into a pair of parallel edges. 
However, adding the triangle within a level-$k$ blob of a level-$k$ network, then adding the triangle would increase the level. Therefore, to prove connectedness of level-$k$ networks under \NNI we use the same idea as for \PR but are more careful to not increase the level.

\begin{theorem}\label{clm:unets:lvlk:connectedness:NNI}
Let $n \geq 3$ and $k \geq 1$.\\
Then, the space $\ulvlknets$ is connected under \NNI with unbounded diameter.
\begin{proof}
Let $N \in \ulvlknets$ and let $T \in \utrees$. 
Like in the proof of \cref{clm:unets:lvlk:connectedness:TBR}, we want to transform $N$ into a network $M$ obtained from $T$ by adding a $k$-burl to the edge incident to leaf $1$.

Let $B$ be a level-$k$ blob of $N$. Assume that $N$ contains another blob $B'$. 
By \cref{clm:unets:TBR:pathDown} there is a \PRP-sequence that can remove $B'$. 
Use \cref{clm:PRMtoNNIM} to substitute this sequence with an \NNI-sequence that reduces $B'$ to a level-1 blob.
Note that this can be done locally within blob $B'$ and its incident edges. 
Therefore, this process does not increase the level of a network along this sequence.
If $B'$ is now a cycle of size at least three, then we can shrink it to a triangle, if necessary, and remove it with an \NNIM.
If $B'$ is a pair of parallel edges and one of its vertices is incident to a degree three vertex $v$ that is not part of a level-$k$ blob, then use an \NNIZ to increase the size of $B'$ into a triangle by including $v$ or merge it with the blob containing $v$.
Next, either remove the resulting triangle, or repeat the process above to remove the new blob.
Otherwise, ignore $B'$ for now and continue with another blob of the current network that is neither $B'$ nor $B$.
When this process terminates, we arrive at a network that has only blob $B$, and, potentially, pairs of parallel edges that are incident to both $B$ and a leaf. That is the case since a pair of parallel edges incident to a degree three vertex not in $B$ could be removed with an \NNIZ and an \NNIM.

If the edge incident to leaf $1$ contains a pair of parallel edges or is incident to a degree three vertex not in $B$, then use $k-1$ \NNIP and \NNIZ (or $k$ in the latter case) to create a $k$-burl next to leaf $1$.
Otherwise, if $B$ is incident to three or more cut-edges, then one of them is not incident to leaf $1$ and can be moved to the edge incident to leaf $1$ with an \NNIZ-sequence. If $B$ is incident to two or fewer cut-edges, there is a vertex incident to three cut edges (since $n \geq 3$) and one of them can be moved to the edge incident to leaf $1$ with an \NNIZ-sequence. Then apply the first case again to create a $k$-burl. Finally, remove $B$ and any remaining pair of parallel edges.
This gives a network $M'$ which consists of a tree $T'$ with a $k$-burl at leaf $1$. There is an \NNIZ-sequence from $T'$ to $T$, which is easily converted into a sequence from $M'$ to $M$. Lastly, note that the diameter is unbounded because for each $r\geq 0$, there is a level-$k$ network with $r$ reticulations.
\end{proof}
\end{theorem}

\section{Isometric relations between spaces} \label{sec:isometric} 
Recall that a space $\class$ is an isometric subgraph of $\unets$ under a rearrangement operation, say TBR, if the TBR-distance of two networks in $\class$ is the same as their TBR-distance in $\unets$. In this section, we investigate this question for $\utrees$ under \TBR, and for tree-based networks and level-k networks under \TBR and \PR.

We start with $\utrees$. The proof of the following theorem follows the proof by Bordewich \etal~\cite[Proposition 7.1]{BLS17} for their equivalent statement for SNPR on rooted phylogenetic trees and networks closely.

\begin{theorem} \label{clm:unets:TBR:treesIsometric}
The space $\utrees$ is an isometric subgraph of $\unets$ under \TBR. Moreover, every shortest \TBR-sequence from $T \in \utrees$ to $T' \in \utrees$ only uses \TBRZ.
  \begin{proof} 
  Let $\distT$ and $\distN$ be the \TBR-distance in $\utrees$ and $\unets$ respectively. To prove the statement, it suffices to show that $\distT(T, T') = \distN(T, T')$ for every pair $T, T' \in \utrees$. Note that $\distT(T, T') \geq \distN(T, T')$ holds by definition. To prove the converse, let $\sigma = (T = N_0, N_1, \ldots, N_k = T')$ be a shortest \TBR-sequence from $T$ to $T'$. Consider the following colouring of the edges of each $N_i$, for $i \in \set{0, \ldots, k}$. Colour all edges of $T = N_0$ blue. For $i \in \set{1, \ldots, k}$ preserve the colouring of $N_{i-1}$ to a colouring of $N_i$ for all edges except those affected by the \TBR. In particular, an edge that gets added or moved is coloured red, an edge resulting from a vertex suppression is coloured blue if the two merged edges were blue and red otherwise, and the edges resulting from an edge subdivision are coloured like the subdivided edge.
  
  Let $F_i$ be the graph obtained from $N_i$ by removing all red edges. We claim that $F_i$ is a forest with at most $k + 1$ components. Since $F_0  = T$, the statement holds for $i = 0$. If $N_i$ is obtained from $N_{i-1}$ by a \TBRP, then $F_i = F_{i-1}$. If $N_i$ is obtained from $N_{i-1}$ by a \TBRZ or \TBRM, then at most one component gets split. Note that $F_k$ is a so-called agreement forest for $T$ and $T'$ and thus $\distT(T, T') \leq k = \distN(T, T')$ by Theorem~2.13 by Allen and Steel~\cite{AS01}. Furthermore, if $\sigma$ would use a \TBRP, then the forest $F_k$ would contain at most $k$ components. However, then $\distT(T, T') < k$; a contradiction.
  \end{proof}
\end{theorem}

\begin{figure}[htb]
  \centering
  \includegraphics{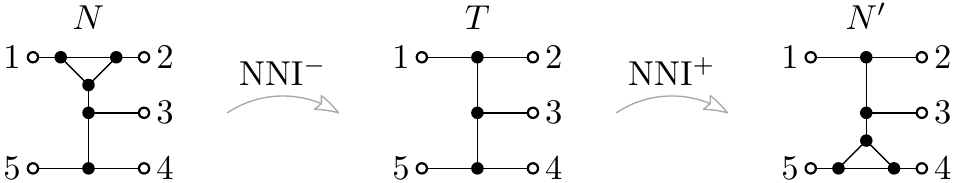}
  \caption{An \NNI-sequence from $N$ to $N'$ using an \NNIP that adds $f$, an \NNIZ that moves $e$, and an \NNIM that removes $e'$. A shortest \NNIZ-sequence from $N$ to $N'$ has length three.}
  \label{fig:unets:NNI:tierNonIsometric} 
\end{figure}

Francis \etal~\cite{FHMW17} gave the example in \cref{fig:unets:NNI:tierNonIsometric} to show that the tiers $\unetsr$ for $n \geq 5$ and $r > 0$ are not isometric subgraphs of $\unets$ under \NNI. Their question of whether tier zero, $\utrees$, is an isometric subgraph of $\unets$ under \NNI remains open.

\begin{lemma} \label{clm:unets:NNI:tierNonIsometric}
Let $n \geq 5$ and $r \geq 0$. Then the space $\unetsr$ is not an isometric subgraph of $\unets$ under \NNI.
\end{lemma}

\begin{figure}[htb]
  \centering
  \includegraphics{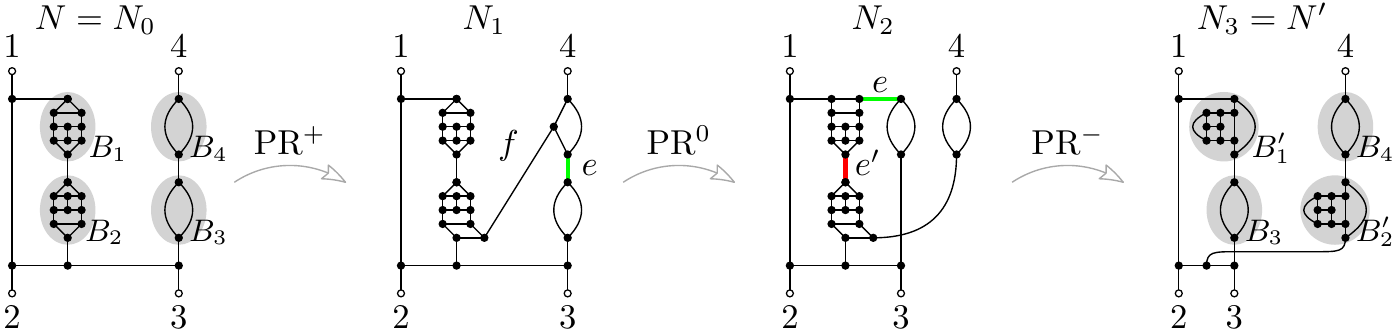}
  \caption{A length three \PR-sequence from $N$ to $N'$ that uses a \PRP, which adds $f$, a \PRZ, which moves $e$, and a \PRM, which removes $e'$. 
  A \PRZ-sequence from $N$ to $N'$ has length at least four.}
  \label{fig:unets:PR:nonIsometric}
\end{figure}

\begin{lemma} \label{clm:unets:PR:tierNonIsometric}
For $n=4$ and $r=13$ the space $\unetsr$ is not an isometric subgraph of $\unets$ under \PR.
  \begin{proof}
  For the networks $N$ and $N'$ in $\unetsr$ shown in \cref{fig:unets:PR:nonIsometric} there is a length three \PR-sequence that traverses tier $r+1$, for example, like the depicted sequence $\sigma = (N = N_0, N_1, N_2,$ $N_3 = N')$. To prove the statement we show that every \PRZ-sequence from $N$ to $N'$ has length at least four.
  
  The networks $N$ and $N'$ contain the highlighted (sub)blobs $B_1$, $B_2$, (resp. $B_1'$ and $B_2'$), $B_3$, and $B_4$.  Observe that the edges between $B_1$ and $B_2$ and between $B_3$ and $B_4$ may only be pruned from a blob by a \PRZ if they get regrafted to the same blob again. Otherwise the resulting network is improper. Note that to derive $B_1'$ from $B_1$ an edge has to be regrafted to the ``top'' of $B_1$ and the edge to $B_2$ has to be pruned. By the first observation, combining these into one \PRZ cannot build the connection to $B_3$. The same applies for the transformation of $B_2$ into $B_2'$ and its connection to $B_4$. Therefore, we either need four \PRZ to derive $B_1'$ and $B_2'$ or two \PRZ plus two \PRZ to build the connections to $B_3$ and $B_4$. In conclusion, at least four \PRZ are required to transform $N$ into $N'$, which concludes this proof.
  \end{proof}
\end{lemma}

By replacing a leaf with a tree, and adding more pairs of parallel edges to edge leading to $4$, this example can be made to work for $n\geq 4$ and $r\geq 13$.

\begin{theorem} \label{clm:unets:tbased:nonIsometric}
For $n \geq 6$ the space $\utbasednets$ is not an isometric subgraph of $\unets$ under \TBR and \PR.
\begin{proof}
  Let $N$ be the network in \cref{fig:unets:tbased:nonIsometric}. Let $N'$ be the network derived from $N$ by swapping the labels $1$ and $2$.
  Note that $\dTBR(N, N') = \dPR(N, N') = 2$, since, from $N$ to $N'$, we can move leaf 2 next to leaf 1 and then move leaf 1 to where leaf 2 was. 
  However, then the network in the middle is not tree-based, since the blob derived from the Petersen graph has no Hamiltonian path if the two pendent edges of the blob are next to each other~\cite{FHM18}.
  We claim that there is no other length two \TBR-sequence from $N$ to $N'$.
  For this proof we call a blob derived from the Petersen graph a Petersen blob.

  \begin{figure}[htb]
    \centering
    \includegraphics{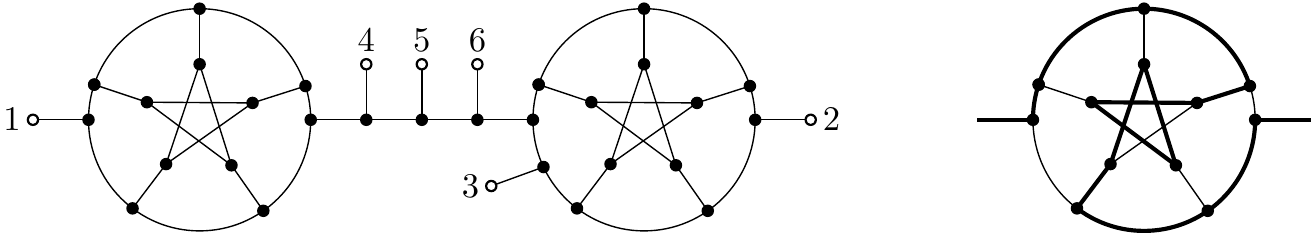}
    \caption{A tree-based network on the left and a Hamiltonian path through a blob derived from the Petersen graph on the right.}
    \label{fig:unets:tbased:nonIsometric}
  \end{figure}
  
  First, note that the \TBRZ-sequence of $N$ and $N'$ is at least two and there is thus no \TBR-sequence that consists of a \TBRM and a \TBRP. Otherwise, these two operations could be merged into a single \TBRZ by \cref{clm:unets:TBR:PMtoZ}. Note that we can only move leaf 1 or 2 by pruning an incident edge if we do not affect the split 1 versus 2, 3 or break the tree-based property. Therefore, they either have to be swapped using edges of the Petersen blobs or the $(4, 5, 6)$-chain has to be reversed and leaf 3 moved to the other Petersen blob.
  However, it is straightforward to check that neither can be done with two \TBRZ. In particular, we can look at what edge the first \TBRZ might move and then check whether a second \TBRZ can arrive at $N'$. If the first \TBRZ breaks a Petersen blob, the problem is that the second \TBRZ has to restore it. We then find that this does not allows us to make the initially planned changes to arrive at $N'$. On the other hand, if we avoid breaking the Petersen blob and reverse the $(4, 5, 6)$-chain, then leaf 3 is still on the wrong side; and if we move leaf 3 to the other Petersen blob, then not enough \TBRZ moves remain to reverse the chain.
  
  Since there is no other length two \TBRZ-sequence there is also no other length two \PR-sequence.
\end{proof}
\end{theorem}

\begin{theorem} \label{clm:unets:lvlk:nonIsometric}
For $n\geq 5$ and large enough $k$, the space $\ulvlknets$ is not an isometric subgraph of $\unets$ under \TBR and \PR.
\begin{proof}
For even $k$, the networks $N$ and $N'$ in \cref{fig:unets:lvlk:nonIsometric} have \TBR- and \PR-distance two via the network $M$. However, note that in $M$ the blobs of size $\frac{k}{2} + 1$ a $\frac{k}{2}$ are merged into a blob of size $k + 1$. Therefore, $M$ is not a level-$k$ network.
We claim that there is no \TBR- or \PR-sequence of length two that does not go through a level-$(k+1)$ network like $M$. An example for odd $k$ can be derived from this.

\begin{figure}[htb]
  \centering
  \includegraphics{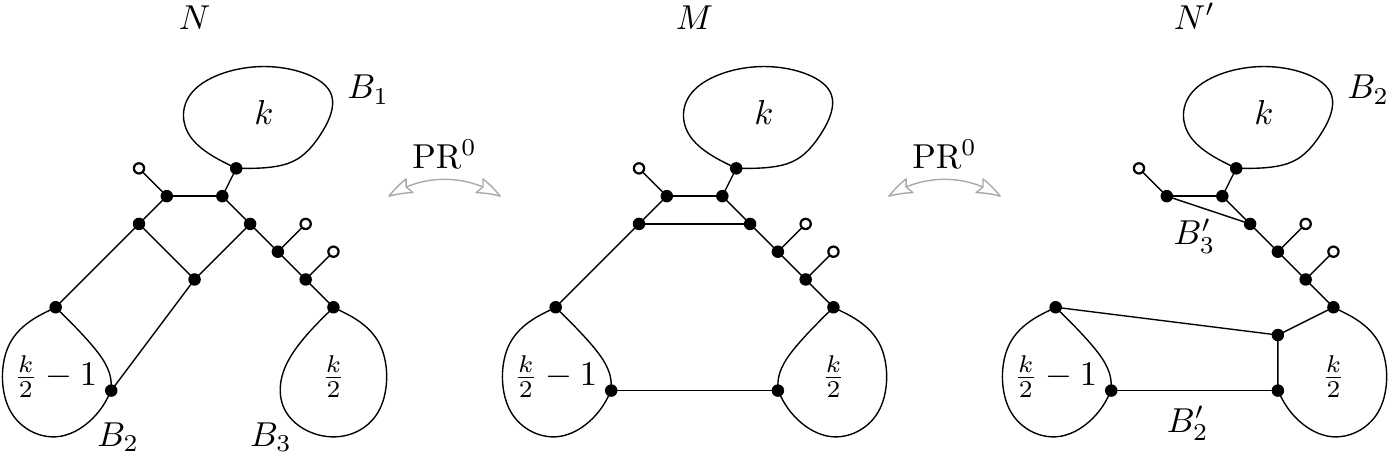} 
  \caption{For even $k$, a \PRZ-sequence from a level-$k$ network $N$ to a level-$k$ network $N'$ (hidden reticulations of the blob-parts given inside, at least two leaves ommited: in $B_1$ and in $B_3$). However, the network $M$ in the middle is a level-$(k+1)$ but not a level-$k$ network.}
  \label{fig:unets:lvlk:nonIsometric}
\end{figure}

It is easy to see that the \TBR-distance of $N$ and $N'$ is at least two and there is thus no \TBR-sequence that consists of a \TBRM and a \TBRP. Otherwise, these two operations could be merged into a single \TBRZ by \cref{clm:unets:TBR:PMtoZ}. We thus have to prove that there is no length two \TBRZ-sequence from $N$ to $N'$ that avoids a level-$(k+1)$ network. Note that it requires two \TBRZ (or \PRZ) to connect $B_2$ and $B_3$ into $B_2'$. Similarly, it requires either two prunings from the upper five-cycle of $B_2$ to obtain the triangle $B_3'$ or one pruning within that cycle. However, in the latter option this would not contribute to connecting $B_2$ and $B_3$ and hence overall at least three operations would be needed. Therefore we have to combine the two operations necessary to create $B_2'$ and to create $B_3'$, which however gives us a sequence like the one shown in \cref{fig:unets:lvlk:nonIsometric}.
\end{proof}
\end{theorem}

Note that the results of this section that show that the spaces of tree-based networks and level-$k$ networks are not isometric subgraphs of the space of all networks also hold if we restrict these spaces to a particular tier $r$ (for large enough $r$).

\section{Computational complexity} \label{sec:complexity} 
In this section, we consider the computational complexity of computing the \TBR-distance and the \PR-distance. First, we recall the known results on phylogenetic trees.

\begin{theorem}[\cite{DGHJLTZ97,HDRB08,AS01}] 
\label{clm:unets:distancesNPhard}
Computing the distance of two trees in $\utrees$ is NP-hard for the \NNI-distance, the \SPR-distance, and the \TBR-distance.
\end{theorem}

In \cref{clm:unets:TBR:treesIsometric}, we have shown that $\utrees$ is an isometric subgraph of $\unets$ under \TBR. Hence, with \cref{clm:unets:distancesNPhard}, we get the following corollary.
\begin{corollary} \label{clm:unets:TBR:NP}
Computing the \TBR-distance of two arbitrary networks in $\unets$ is NP-hard.
\end{corollary}

We can use the same two theorems to prove that computing the \TBR-distance in tiers is also hard.

\begin{theorem} \label{clm:unets:TBR:tierNP}
Computing the \TBR-distance of two arbitrary networks in $\unetsr$ is NP-hard.
  \begin{proof}
    We (linear-time) reduce the NP-hard problem of computing the \TBR-distance of two trees in $\utrees$ to computing the \TBR-distance of two networks in $\unetsxx[n+1,r]$. For this, let $T, T' \in \utrees$. Let $e$ be the edge incident to leaf $n$ of $T$. Obtain $S$ from $T$ by subdividing $e$ with a new vertex $u$ and adding the edge $\{u, v\}$ where $v$ is a new vertex labelled $n+1$. Next, add $r$ handcuffs to the cherry $\set{n, n+1}$ to obtain the network $N \in \unetsxx[n+1,r]$. Analogously obtain $N'$ from $T'$.
    
    The equality $\dTBR(T, T') = \dTBR(N, N')$ follows from \cref{clm:unets:TBR:existingCloseDisplayedTree}, and the fact that networks handcuffed at a cherry display exactly one tree. More precisely, a \TBR-sequence between $T$ and $T'$ induces a \TBR-sequence of the same length between $N$ and $N'$, hence $\dTBR(T, T') \geq \dTBR(N, N')$.
    Conversely, by \cref{clm:unets:TBR:existingCloseDisplayedTree} and the fact that $D(N)=\{T\}$ and $D(N')=\{T'\}$, it follows that $\dTBR(T, T') \leq \dTBR(N, N')$. Since computing the TBR-distance in $\utrees$ is NP-hard, the statement follows.
  \end{proof}
\end{theorem}

To prove that computing the \PR-distance is hard, we use a different reduction. Van Iersel et al. prove that deciding whether a tree is displayed by a (not necessarily proper) phylogenetic network (Unrooted Tree Containment; UTC) is NP-hard \cite{vIKSSB18}. Combining this with \cref{clm:unets:TBR:pathDown}, we arrive at our result.

\begin{theorem}
Computing the \PR-distance of two arbitrary networks in $\unets$ is NP-hard.
  \begin{proof}
    We reduce from UTC to the problem of computing the \PR-distance of two networks in $\unets$. Let $(N,T)$ with $N$ a (not necessarily proper) network and $T\in\utrees$ be an arbitrary instance of UTC. We obtain an instance $(N',T',r')$ of the \PR-distance decision problem as follows: remove all cut-edges of $N$ that do not separate two labelled leaves, and let $N''$ be the connected component containing all the leaves; now, let $N'$ be the proper network obtained from $N''$ by suppressing all degree two nodes. The instance of the \PR-distance decision problem consists of $N'$, $T'=T$, and the reticulation number $r'$ of $N'$. As we can compute in polynomial time whether a cut edge separates two labelled leaves, the reduction is polynomial time. Because a displayed tree uses only cut-edges that separate two labelled leaves, $T$ is displayed by $N$ if and only if it is displayed by $N'$. By \cref{clm:unets:TBR:pathDown}, $T$ is a displayed tree of $N$, if and only if $\dPR(N',T')\leq r$, which concludes the proof.
  \end{proof}
\end{theorem}

Unlike for the hardness proof of \TBR-distance, we cannot readily adapt this proof to the \PR-distance in $\unetsr$. For this purpose, we need to learn more about the structure of \PR-space.

\section{Concluding remarks} 
In this paper, we investigated basic properties of spaces of unrooted phylogenetic networks and their metrics under the rearrangement operations \NNI, \PR, and \TBR.
We have proven connectedness and bounds on diameters for different classes of phylogenetic networks, including networks that display a particular set of trees, tree-based networks, and level-$k$ networks.
Although these parameters have been studied before for classes of rooted phylogenetic network~\cite{BLS17}, this is the first paper that studies these properties for classes of unrooted phylogenetic networks besides the space of all networks. A summary of our results is shown in \cref{tbl:unets:connectedness}.

To see the improvements in diameter bounds, we compare our results to previously found bounds: For the space of phylogenetic trees $\utrees$ it was known that the diameter is asymptotically linearithmic and linear in the size of the trees under \NNI and \SPR/\TBR~\cite{LTZ96,DGH11}, respectively. Here, we have shown that the diameter under \NNI is also asymptotically linearithmic for higher tiers of phylogenetic networks. Whether this also holds in the rooted case is still open. We have further (re)proven the asymptotic linear diameter for \PR and \TBR of these tiers and, in particular, improved the upper bound on the diameter under \TBR to $n - 3 - \floor{\frac{\sqrt{n - 2} - 1}{2}} + r$ from the previously best bound $n + 2r$~\cite{JJEvIS17}. 

\begin{table}[htb]
    \centering
    \begin{tabular}{c|c|c|c}
        class & \NNI & \PR & \TBR \\ \hline
        $\utrees$ & $\Theta(n \log n)$ \cite{LTZ96}  & $\Theta(n)$ \cite{DGH11} & $\Theta(n)$  \cite{DGH11} \\
        $\unetsr$ & $\Theta(m \log m)$ T.~\ref{clm:unets:NNIconnected:tier} 
            & $\Theta(m)$ \cite{FHM18,JJEvIS17} & $\Theta(m)$ T.~\ref{clm:unets:TBR:tierConnected} \\
        $\unets$  & \checkmark \cref{clm:unets:NNI:connected} & \checkmark  \cref{clm:unets:TBR:connected} & \checkmark  \cref{clm:unets:TBR:connected} \\
        $\unets(P)$ & \checkmark \cref{clm:unets:displayingSet:connectivity} & \checkmark \cref{clm:unets:displayingSet:connectivity} & \checkmark \cref{clm:unets:displayingSet:connectivity}\\
        $\utbasednetsr(T)$ & $\Oh(rm)$ \cref{clm:unets:tbased:connectedness} & $\Theta(r)$ \cref{clm:unets:tbased:connectedness} & $\Theta(r)$ \cref{clm:unets:tbased:connectedness} \\
        $\utbasednetsr$    & $\Oh(rm + n \log n)$ T.~\ref{clm:unets:tbased:connectednessOther} & $\Theta(m)$ \cref{clm:unets:tbased:connectednessOther} & $\Theta(m)$ T.~\ref{clm:unets:tbased:connectednessOther} \\
        $\utbasednets(T)$  & \checkmark \cref{clm:unets:tbased:connectednessOther} & \checkmark \cref{clm:unets:tbased:connectednessOther} & \checkmark \cref{clm:unets:tbased:connectednessOther} \\
        $\utbasednets$     & \checkmark \cref{clm:unets:tbased:connectednessOther} & \checkmark \cref{clm:unets:tbased:connectednessOther} & \checkmark \cref{clm:unets:tbased:connectednessOther} \\
        $\ulvlknets$   & \checkmark \cref{clm:unets:lvlk:connectedness:NNI} & \checkmark \cref{clm:unets:lvlk:connectedness:TBR} & \checkmark \cref{clm:unets:lvlk:connectedness:TBR} \\  
    \end{tabular}
    \caption{Connectedness and diameters, if bounded, for the various classes and rearrangement operations. Here $m = n + r$, $P$ is a set of phylogenetic networks, and $T \in \utrees$.}
    \label{tbl:unets:connectedness}
\end{table}

To uncover local structures of network spaces, we looked at properties of shortest sequences of moves between two networks. Here we found that shortest \TBR-sequences between networks in the same tier never traverse lower tiers, and shortest \TBR-sequences between trees also never traverse higher tiers. This implies that $\utrees$ is an isometric subgraph of $\unets$, and that computing the \TBR-distance between two networks in $\unets$ is NP-hard. This answers a question by Francis \etal~\cite{FHMW17}. 
We have attempted to prove similar results for other subspaces and rearrangement moves. However, for higher tiers, we have not been able to prove that shortest \TBR-sequences never traverse higher tiers. To answer this question we may need to utilise agreement graphs such as frequently used for phylogenetic trees~\cite{AS01,BS05} and, more recently, also for rooted phylogenetic networks~\cite{KL19,Kla19}.
Concerning \NNI and \PR we gave counterexamples to prove that higher tiers are not isometric subgraphs of $\unets$. The questions whether $\utrees$ is isometrically embedded in $\unets$ under \PR and \NNI remains open. Answering these questions positively would also provide an answer to the question whether computing the shortest \NNI-distance between two networks is NP-hard, and clues toward proving whether the \PR-distance between two networks in the same tier is NP-hard. Further negative results that we have shown are that the spaces of tree-based networks and level-$k$ are not isometric subgraphs of the space of all phylogenetic networks.

Throughout this paper, we have restricted our attention to proper networks. We could also have chosen to use unrooted networks without the properness condition. This definition, which is mathematically more elegant, is used in most other papers, so it seems to be the obvious choice. However, it is not natural to have cut-edges that do not separate leaves: such networks carry no biological meaning. It is desirable that networks are rootable and thus have an evolutionary interpretation. Unrooted phylogenetic networks are rootable if they have at most one blob with one cut-edge. While using this in the definition of an unrooted phylogenetic network could therefore be sufficient, we go one step further, and ask that there is no such blob. This makes a network rootable at any leaf (i.e., with any taxon as out-group), which gives a stronger biological interpretation and usability.

The fact that our definition of unrooted phylogenetic networks is mathematically more restrictive, means that any positive result we have proven is likely also true when using a less restrictive definition. That is, connectedness for those definitions follows easily by finding sequences to proper networks, like done by Jansen \etal~\cite{JJEvIS17}.
As we may be able to find short sequences for this purpose, the diameter results will likely also still hold. This means that whatever definitions may be used in practice, with minor additional arguments, our results provide the theoretical background necessary to justify local search operations.

\pdfbookmark[1]{Acknowledgments}{Acknowledgments}
\subsection*{Acknowledgements} 
The first author was supported by the Netherlands Organization for Scientific Research (NWO) Vidi grant 639.072.602.
The second author thanks the New Zealand Marsden Fund for their financial support.

\phantomsection
\pdfbookmark[1]{References}{references} 

\providecommand{\etalchar}[1]{$^{#1}$}

\end{document}